\documentclass{article}
\usepackage{graphicx} 
\usepackage[utf8]{inputenc}
\usepackage{amsmath}
\usepackage{amssymb}
\usepackage{amsthm}
\usepackage{dsfont}
\usepackage[a-1b]{pdfx} 
\usepackage{hyperref}
\usepackage{csquotes}
\usepackage{pdflscape}
\usepackage{algorithm}
\usepackage{algpseudocode}
\usepackage[a4paper, total={6in, 8.5in}]{geometry}
\hypersetup{colorlinks,linkcolor={black},citecolor={blue},urlcolor={black}}  
\usepackage[english]{babel}
\usepackage{comment}
\usepackage{bbm}
\usepackage{float}
\usepackage{rotating}
\usepackage{array}
\usepackage{multirow}

\usepackage[round]{natbib}
\usepackage[toc,page]{appendix}
\renewcommand{\arraystretch}{1}
\usepackage{authblk}

\newtheorem{theorem}{Theorem}[section]
\newtheorem{corollary}[theorem]{Corollary}
\newtheorem{lemma}[theorem]{Lemma}
\newtheorem{proposition}[theorem]{Proposition}
\theoremstyle{definition}

\newtheorem{definition}[theorem]{Definition}
\newtheorem{example}[theorem]{Example}

\newtheoremstyle{mystyle}
  {\topsep}
  {\topsep}
  {\normalfont}
  {}
  {\bfseries}
  {.}
  {.5em}
  {}

\theoremstyle{mystyle}
\newtheorem{remark}[theorem]{Remark}

\newcommand{\pr}{\mathbb{P}}
\newcommand{\p}{\mathbb{P}}

\newcommand{\Ev}{\mathbb{E}}
\newcommand{\ind}{\mathds{1}}

\renewcommand{\d}{\mathrm{d}}

\title{Combining exchangeable p-values}
\author[1]{Matteo Gasparin}
\author[2]{Ruodu Wang}
\author[3,4]{Aaditya Ramdas}

\affil[1]{Department of Statistical Sciences, University of Padova}
\affil[2]{Department of Statistics and Actuarial Science, University of Waterloo}
\affil[3]{Department of Statistics and Data Science, Carnegie Mellon University}
\affil[4]{Machine Learning Department, Carnegie Mellon University}
\date{\today}

\begin{document}

\maketitle

\begin{abstract}
The problem of combining p-values is 
an old and fundamental one,
and the classic assumption of independence is often violated or unverifiable in many applications. 
There are many well-known rules that can combine a set of arbitrarily dependent p-values (for the same hypothesis) into a single p-value. 
We show that essentially all these existing rules can be strictly improved when the p-values are exchangeable, or when external randomization is allowed (or both). 
For example, we derive randomized and/or exchangeable improvements of well known rules like ``twice the median'' and ``twice the average'', as well as geometric and harmonic means.
 Exchangeable p-values are often produced one at a time (for example, under repeated tests involving data splitting), and our rules can combine them sequentially as they are produced, stopping when the combined p-values stabilize. Our work also improves rules for combining arbitrarily dependent p-values, since the latter becomes exchangeable if they are presented to the analyst in a random order.
The main technical advance is to show that all existing combination rules can be obtained by calibrating the p-values to e-values (using an $\alpha$-dependent calibrator), averaging those e-values, converting to a level-$\alpha$ test using Markov's inequality, and finally obtaining p-values by combining this family of tests; the improvements are delivered via recent randomized and exchangeable variants of Markov's inequality.
\end{abstract}

\tableofcontents

\section{Introduction}\label{sec:intro}
The combination of p-values represents a fundamental task frequently encountered in statistical inference and its applications in the natural sciences. Within the realm of multiple testing, for instance, the focus lies in testing whether all individual null hypotheses are simultaneously true. This particular challenge, often referred to as global null testing, can be addressed by merging multiple p-values into a single p-value. Potential solutions, assuming the p-values are statistically independent, are provided in \cite{fisher1928,pearson1934,simes1986}, with the latter also working under a certain notion of positive dependence~\citep{sarkar1998, benjamiyekuteli2001}. See \cite{owen2009} for a review of the methods. Recently, harmonic mean p-values~\citep{wilson2019} and methods under negative dependence~\citep{chi2022multiple} have been developed in the context of multiple testing.

The assumption of independence (or positive/negative dependence) is often violated in many real-world applications (see e.g., Section 4.2 of \citealp{efron2012large}), and in certain scenarios, it may be preferable not to impose unverifiable conditions on the joint distribution of the p-values, beyond the minimum necessary assumption that each individual p-value is indeed (marginally) valid. Several methods are available for combining p-values with arbitrary dependencies; notably, the Bonferroni method is widely used, which involves multiplying the minimum of the p-values by the number of tests conducted. Other methods have been proposed in the literature, some based on order statistics \citep{ruger1978,morgenstern1980,hommel1983}, while others rely on their arithmetic mean and other variants \citep{ruschendorf1982,vovk2020, vovk2022}.  Two prominent examples are that both 2 times the median of the p-values, and 2 times the average of the p-values, are valid combination rules, and the multiplicative factor of 2 cannot be reduced. Inevitably, these methods satisfying  validity under arbitrary dependence come with a price to pay in terms of statistical power. Our work will improve all of these rules under the weaker assumption of validity under exchangeability of p-values.

To elaborate, the main objective of our work is to obtain simple new valid merging methods under the assumption of exchangeability of the input p-values, which are more powerful than methods that assume arbitrary dependence. Implied by the relative strength of the dependence assumptions, the new methods will be incomparable to methods assuming negative or positive dependence, and less powerful than methods assuming independence. However, specialized methods for handling exchangeable dependence are quite practically relevant. Such dependence is encountered, for example, in statistical testing via sample splitting, as we now elaborate. 

There are at least two different reasons for which sample splitting is used: the first is to relax the assumptions needed to obtain theoretical guarantees, while the second is to reduce computational costs. Some examples of such procedures are \cite{cox1975,wasserman2009,banejeree2019,wasserman2020pnas,shafer2008,shekhar2022permutation,shekhar2023permutation,Kim2024}. The drawback of these methods based on sample-splitting is that the obtained p-values are affected by the randomness of the split. \cite{meinshausen2009} called this phenomenon as a \emph{p-value lottery}. So one may instead repeat the same sample splitting procedure several times to obtain multiple p-values, which are exchangeable by design of the procedure.

One can of course combine such p-values by using the earlier mentioned rules (like twice the average) for arbitrarily dependent p-values. But we hope to do better by exploiting their exchangeability. In a paper that seemingly dampens that hope, \cite{choi2023} showed that in the aforementioned rule of ``twice the average'', the constant factor of 2 cannot be improved even under exchangeability. However, their result does not imply that ``twice the average'' cannot be improved; it simply states that any such improvement cannot proceed by attempting to lower the constant of 2. Indeed, our paper will improve on this well known rule (and many others), but it proceeds differently, not by lowering the constant. Instead, it calculates the twice the average of the first $k$ p-values, and takes a minimum over all $k$ (see Table~\ref{tab:1}).

We also show that no symmetric rule for merging arbitrarily dependent p-values (like the ones mentioned earlier) can be improved under exchangeability. To achieve any improvement, we must consider asymmetric rules that process the p-values in a particular order. This might initially appear paradoxical given the exchangeability of the p-values, but it is easily sorted out. 
In many practical settings, these exchangeable p-values can be generated one by one by repeating the same randomized procedure many times, generating a stream of p-values. 
In this case, our combination rules would simply process these p-values in the order that they are generated. This seems quite appropriate, and the main advantage of doing this is increased power over processing them symmetrically as a batch.
In fact, as we discuss, we do not need to fix the number of p-values ahead of time, they can just be processed online, yielding a p-value whenever this procedure is stopped. This makes our merging rules particularly simple and practical.

The problem of combining such p-values from repeated sample splitting has been studied by other authors, such as \cite{diciccio2020}, but our combination rules are more powerful, 
and also more general and systematic because they apply more broadly. The same problem has also been studied in \cite{guo2023}, where the authors propose a combination method based on subsampling to merge test statistics based on different random splits. Unlike our nonasymptotic guarantees that work directly with the p-values and are cheap to compute, their work provides only asymptotic guarantees under certain additional assumptions (like an asymptotic pivotal null distribution for their test statistics) and they require access to the full dataset on which to perform expensive subsampling-based recomputations. However, when their additional assumptions hold, their procedure can be expected to be more powerful than ours because it estimates and exploits the joint distribution of the p-values. With a similar aim, \cite{ritzwoller2023} investigates a method to enhance the reproducibility of statistical results obtained through sample-splitting. Specifically, their algorithm sequentially aggregates statistics across multiple sample splits until the variability induced by the different splits falls below a defined threshold.

It is perhaps interesting that all the aforementioned methods for merging under arbitrary dependence or under exchangeability are actually inadmissible when randomization is permitted. Randomization in the context of hypothesis testing is not new and is used, for example, in discrete tests; some examples are Fisher's exact test \citep{fisher1928} or the randomized test for a binomial proportion proposed in \cite{stevens1950}. In our paper, we will see how the introduction of a simple external randomization (an independent uniform random variable and/or uniform permutation) can improve the existing merging rules for arbitrary dependence, as well as our new rules for exchangeable merging. 

In terms of technical aspects, one of our main contributions is to point out explicitly how existing merging rules for arbitrary dependence are actually recovered in a unified manner: by transforming the p-values into e-values \citep{vovk2021evals,wasserman2020pnas,grunwald_safe_2019} using different ``calibrators'', averaging the resulting e-values and finally applying Markov's inequality. This connection is particularly important, because then the improvements under exchangeability, or by randomization, are then achieved by invoking the recent ``exchangeable Markov inequality'' and ``uniformly randomized Markov inequality''~\citep{ramdas2023}.

\paragraph{Paper outline and peek at results.}
The rest of this paper is organized as follows. In Section~\ref{sec:notation}, we introduce the notation and tools necessary for the paper. In Section~\ref{sec:general_res}, the main results are presented in a general way, focusing on two distinct aspects: the first part addresses the case of exchangeable p-values, while the second part introduces novel findings under the assumption of arbitrarily dependent p-values when randomization is allowed. Subsequent sections investigate the implications of these results across various p-merging functions commonly found in the literature. Specifically, Section~\ref{sec:ruger} and Section~\ref{sec:hommel} delve into the combination proposed by~\cite{ruger1978} and~\cite{hommel1983}, respectively. Section~\ref{sec:avg} examines the case of arithmetic mean. The following two sections address two additional scenarios within the family of generalized means: namely, the harmonic mean (Section~\ref{sec:harm_avg}) and the geometric mean (Section~\ref{sec:geo_avg}).
Section~\ref{sec:simul} presents some simulation results, before we conclude in Section~\ref{sec:summary}. All proofs are provided in Appendix \ref{sec:ap_proof}.

Before proceeding with the paper, Table~\ref{tab:1} presents some notable combination rules introduced in the literature and their corresponding exchangeable and randomized versions introduced in the following sections. 
These results are first derived in a general form and then discussed case-by-case. Some of the rules in the table are not admissible, as will be explained in the following.

\renewcommand{\arraystretch}{3}
\begin{table}[t]
    \centering
    \begin{tabular}{l|c|c|c}
    Combination rule  & {\begin{minipage}{0.22\textwidth}\centering 
 Arbitrary dependence \\ (known)
    \end{minipage}}
    & {\begin{minipage}{0.18\textwidth}\centering 
 Exchangeability\\ (new)
    \end{minipage}}
    & 
    {\begin{minipage}{0.22\textwidth}\centering 
 Arbitrary  dependence,\\  randomized (new)
    \end{minipage}}
    \\ \hline
    R\"uger combination & $\frac{K}{k} p_{(k)}$ & $\frac{K}{k} \bigwedge_{m=1}^K p^m_{(\lambda_m)}$ & $\frac{K}{k} p_{(\lceil Uk \rceil)}$ 
     \\ Arithmetic mean &  $2 A(\mathbf{p})$ & $2 \bigwedge_{m=1}^K A(\mathbf{p}_m)$ &   $\frac{2}{2-U} A(\mathbf{p})$
    \\ Geometric mean &  $e G(\mathbf{p})$   &   $e \bigwedge_{m=1}^K G(\mathbf{p}_m)$ & $e^U G(\mathbf{p})$ 
    \\ Harmonic mean & $T_K^+ H(\mathbf{p})$   & $T_K^+ \bigwedge_{m=1}^K H(\mathbf{p}_m)$   &  $(T_KU + 1)H(\mathbf{p})$ \\
    \end{tabular}

    \label{tab:1}
    \caption{Some combination rules for arbitrarily dependent p-values documented in literature, along with their exchangeable and randomized improvements. If randomization is permitted, one can also improve the existing rules for combining arbitrarily dependent p-values by using the exchangeable combination rule applied to a random permutation of the p-values (this is not presented as a separate column). Here, $\mathbf{p}=(p_1, \dots, p_K)$ denotes the vector of p-values, and $\mathbf{p}_m$ represents the vector containing the first $m$ values of $\mathbf{p}$. In the table, $p_{(k)}$ is the $k$-th smallest value of $\mathbf{p}$, while $p_{(\lambda_m)}^m$ is the $\lambda_m=\lceil m\frac{k}{K} \rceil$ ordered value of $\mathbf{p}_m$. 
    The random variable $U$ is uniformly distributed in the interval $[0,1]$.
    Additionally, $A, G$ and $H$ respectively denote the arithmetic mean, the geometric mean, and the harmonic mean. The values $T_K$ and $T_K^+$ are given by $T_K=\log K + \log \log K +1$ and $T_K^+=T_K+1$, for $K \geq 2$.}
\end{table}

\section{Problem setup and notation}\label{sec:notation}
Without loss of generality, let $(\Omega,\mathcal F, \pr)$ be an atomless probability space\footnote{A probability space is atomless if there exists a random variable  on this space that is uniformly distributed on $[0,1]$.}, and this is implicitly assumed in almost all papers in statistics; see \citet[Appendix D]{vovk2021evals} for related results and discussions. Let $\mathcal U$ be the set of all uniform random variables on $[0,1]$ under $\pr$.
In the following, $K \geq 2$ is an integer. We use the shorthand notation $x \vee y = \max(x,y)$, $x \wedge y=\min(x,y)$,
 $\bigvee_{k=1}^K x_k = \max\{x_1,\dots,x_K\}$, and  $\bigwedge_{k=1}^K x_k= \min\{x_1,\dots,x_K\}$.

A \emph{p-variable}  for testing $\mathbb P$  is a random variable $P: \Omega \to [0, \infty)$ satisfying 
\begin{equation*}
    \pr(P \leq \alpha) \leq \alpha, 
\end{equation*}
for all $\alpha \in (0,1)$. Typically, of course, $P$ will only take values in $[0,1]$, but nothing is lost by allowing the larger range above.
For all results on validity of the methods in this paper, it suffices to consider p-variables in $ \mathcal U$, i.e., $\mathbb P(P\le \alpha)=\alpha$ for each $\alpha \in(0,1)$. 

An \emph{e-variable}  for testing $\mathbb P$ is a non-negative extended random variable $E:\Omega\to[0,\infty]$ with $\Ev_\pr[E]\leq 1$.
A \emph{calibrator} is a decreasing function $f:[0,\infty)\to[0,\infty]$
satisfying $ f=0$ on $(1,\infty)$ and $\int_0^1 f(p)\d p\leq 1$.
Essentially, a calibrator transforms any p-variable to an e-variable.
It is \emph{admissible} if it is upper semicontinuous, $f(0)=\infty$, and $\int_0^1 f(p) \d p = 1$. Equivalently, a calibrator is admissible if it is not strictly dominated, in a natural sense, by any other calibrator~\citep[Proposition~2.1 and Proposition~2.2 in][]{vovk2021evals}.
We fix $\mathbb P$ throughout, and omit ``for testing $\mathbb P$" when discussing p-variables and e-variables; we do not distinguish them from the commonly used terms ``p-values" and ``e-values",    and this should create no confusion.

Our starting point is a collection of $K$ p-variables $\mathbf{P}=(P_1, \dots, P_K)$ and we denote their observed (realized) values by $\mathbf{p}=(p_1, \dots, p_K)$. Borrowing terminology from~\cite{vovk2022} and~\cite{vovk2020}, we have that a \emph{p-merging function} is an increasing Borel function $F:[0,\infty)^{K+1} \to [0,\infty)$ such that $\pr(F(\mathbf{P}) \leq \alpha) \leq \alpha$ whenever $P_1, \dots, P_K$ are p-variables. In other words, the function $F$, starting from $K$ p-values, returns a p-value. A p-merging function is \emph{symmetric} if $F(\mathbf{p})$ is invariant under any permutation of $\mathbf{p}$, and it is \emph{homogeneous} if $F(\gamma \mathbf{p})=\gamma F(\mathbf{p})$ for all $\mathbf{p}$ with $F(\mathbf{p})\leq 1$ and $\gamma \in (0, 1]$. 
The class of homogeneous p-merging functions encompasses the \emph{O-family} based on quantiles introduced in \cite{ruger1978}, the Hommel's combination and the \emph{M-family} introduced in \cite{vovk2020}. We now introduce the notion of \emph{domination} in the context of p-merging functions.
\begin{definition}
    A function $F$ dominates (interpreted as better being smaller) a function $G$ if 
\begin{equation*}
F(\mathbf{p}) \leq G(\mathbf{p}), \quad \mathrm{for~all~}\mathbf{p},
\end{equation*}
and the domination is strict if $F(\mathbf{p}) < G(\mathbf{p})$, for at least one $\mathbf{p}$. 
A p-merging function $F$ is \emph{admissible} if it is not strictly dominated by any other p-merging function.
\end{definition}

Although we defined p-values and e-values for a single probability measure $\mathbb P$, all results hold for composite hypotheses. More precisely, if $\mathbf P$ is a vector of p-variables for a composite hypothesis and $F$ is a p-merging function,
then 
$F(\mathbf P)$ is a p-variable for the same composite hypothesis. See \cite{vovk2021evals} for precise definitions and related discussions.

For any function $F:[0,\infty)^K\to [0,\infty)$ and $\alpha \in (0,1)$, let its rejection region at level $\alpha$ be given by 
\begin{equation*}
  R_\alpha(F) := \left\{\mathbf p \in [0,\infty)^K: F(\mathbf p)\leq \alpha \right\}.
\end{equation*}
For any homogeneous $F$,  $R_\alpha(F)$ for  $\alpha\in(0,1)$
takes the form $R_\alpha(F)=\alpha A$ for some $A\subseteq [0,\infty)^K$,
where $\alpha A$ means the set $\{\alpha\, \mathbf x:\mathbf x\in A\}$.

Conversely, any increasing collection of Borel lower sets
$\{R_\alpha \subseteq [0,\infty)^K: \alpha \in (0,1)\}$
determines an increasing Borel function $F: [0,\infty)^K\to[0,1]$ by the equation
\begin{equation*}
  F(\mathbf p) = \inf\{\alpha\in(0,1): \mathbf p\in R_{\alpha}\},
\end{equation*}
with the convention $\inf\varnothing = 1$ (throughout).
It is immediate to see that $F$ is a p-merging function if and only if
$\pr(\mathbf P \in R_\alpha) \le \alpha$ for all $\alpha\in (0,1)$ and $\mathbf P\in \mathcal U^K$, where $\mathcal U^K$ is the set of all $K$-dimensional random vectors with components in $\mathcal U$. 

Below, $\Delta_K$ is the standard $K$-simplex. 
Every admissible homogeneous p-merging function possesses a dual formulation expressed in terms of calibrators, as summarized below.

\begin{theorem}[\cite{vovk2022}; Theorem 5.1]\label{th:vovkwang}
  For any admissible homogeneous p-merging function $F$,
  there exist $(\lambda_1,\dots,\lambda_K)\in\Delta_K$ and admissible calibrators $f_1,\dots,f_K$
  such that
  \begin{equation}\label{eq:calibrator}
    R_\alpha(F)
    = 
    \left\{
      \mathbf p \in [0,\infty)^K:
      \sum_{k=1}^K \lambda_k f_k\left(\frac{p_k}{\alpha}\right) \ge 1
    \right\}
    \qquad
    \mathrm{for~each~\alpha \in (0,1)}.
  \end{equation}
  Conversely, for any $(\lambda_1,\dots,\lambda_K)\in\Delta_K$ and calibrators $f_1,\dots,f_K$,
  \eqref{eq:calibrator} determines a homogeneous p-merging function.
\end{theorem}
We will exploit this dual form to implement our ``randomized" and ``exchangeable'' techniques, generating p-merging functions that consistently give smaller p-values (usually strictly) than those produced by their \emph{original} counterparts. From~\eqref{eq:calibrator}, it is worth noting that $\sum_{k=1}^K \lambda_k f_k(P_k)$ is an e-value. We now present a useful lemma.
\begin{lemma}\label{lemma:eval}
    Let $f_1, \dots, f_K$ be $K$ calibrators and $\mathbf{P}\in\mathcal{U}^K$. Then, for any $(\lambda_1, \dots, \lambda_K) \in \Delta_K$ and $\alpha \in (0,1]$, 
    \begin{equation}\label{eq:sum_of_cals}
        \frac{1}{\alpha} \sum_{k=1}^K \lambda_k f_k\left(\frac{P_k}{\alpha}\right)
    \end{equation}
    is an e-variable. If $f_1, \dots, f_K$ are admissible calibrators, then 
    \(
    \Ev\left[ \frac{1}{\alpha} \sum_{k=1}^K \lambda_k f_k\left(\frac{P_k}{\alpha}\right) \right] = 1.
    \)
\end{lemma}
Lemma \ref{lemma:eval} used the fact that a calibrator takes value $0$ on $(1,\infty)$, which is not restrictive because the relevant range of p-values is $[0,1]$. This condition is   assumed when defining calibrators in \cite{vovk2022}.

In particular, choosing $\lambda_1=1$ and $\lambda_k=0$, for $k\geq2$, we have that $(1/\alpha)f_1(P_1/\alpha)$ is an e-value, for all $\alpha \in (0,1]$. 

Before introducing our results, we present some inequalities introduced in~\cite{ramdas2023} that will be fundamental throughout the subsequent discussion. The following inequalities can be viewed as an extension of Markov's inequality. 
\begin{theorem}[Exchangeable Markov Inequality]\label{th:EMI}
Let $X_1, X_2, \dots$ form an exchangeable sequence of non-negative and integrable random variables. Then, for any $a>0$,
\[
\pr\left(\exists k \geq 1: \frac{1}{k} \sum_{i=1}^k X_i \geq \frac{1}{a}\right) \leq a \Ev[X_1].
\]
In addition, let $X_1, \dots, X_K$ be exchangeable, non-negative and integrable random variables. Then, for any $a>0$,
\[
\pr\left(\exists k \leq K: \frac{1}{k} \sum_{i=1}^k X_i \geq \frac{1}{a}\right) \leq a \Ev[X_1].
\]
\end{theorem}

The second inequality is based on an external randomization of the threshold of Markov's inequality. 

\begin{theorem}[Uniformly-randomized Markov Inequality]\label{th:UMI} 
Let $X$ be a non-negative random variable independent of $U \in \mathcal U$. Then, for any $a>0$,
\[
\pr(X \geq U/a) \leq a \Ev[X].
\]  
\end{theorem}
The third inequality combines the previous two theorems in the following way:
\begin{theorem}[Exchangeable and Uniformly-randomized Markov Inequality]\label{th:EUMI}
Let $X_1, \dots, X_K$ be a set of exchangeable and non-negative random variables independent of $U \in \mathcal{U}$. Then, for any $a > 0$, 
\[
\pr\left( X_1 \geq U/a \mathrm{~or~} \exists k \leq K: \frac{1}{k} \sum_{i=1}^k X_i \geq \frac{1}{a} \right) \leq a \, \Ev[X_1].
\]
\end{theorem}

These results will be used in the next section as technical tools to derive new combination rules in different situations.

\section{New results on merging p-values}\label{sec:general_res}

This section introduces our main results stated in abstract and general terms, which we instantiate in special cases (like arithmetic or geometric mean) in the sections that follow.

\subsection{Exchangeable p-values}\label{sec:exch_gen}
Assuming iid~p-values is often overly stringent in numerous practical applications. A more pragmatic and less restrictive assumption is exchangeability of p-values, indicating that the distribution of the p-values is unchanged under a permutation of the indices. Formally, 
\[
(P_1, \dots, P_K) \stackrel{d}{=} (P_{\sigma(1)}, \dots, P_{\sigma(K)}),
\]
where $\stackrel{d}{=}$ represents equality in distribution while $ \sigma:\{ 1,\dots,K\} \to \{1,\dots,K\}$ is any permutation of the indices. As discussed in Section~\ref{sec:intro}, this situation is encountered in statistical testing using repeated 
sample splitting (repeated $K$ times on the same data in an identical fashion). In this section, we assume that the sequence of p-variables $\mathbf{P}=(P_1, \dots, P_K)$ is exchangeable and takes values in $[0,1]^K$. 
\begin{remark}\label{rem:31}
The reader may note that exchangeability can be induced by processing the (potentially non-exchangeable) sequence of p-values $(P_1, \dots, P_K)$ in a uniformly random order. As a consequence, this implies that if randomization is allowed, it is always possible to satisfy the exchangeability assumption even if the starting sequence has an arbitrary dependence.
Stated alternatively, exchangeable combination rules can be safely applied  to arbitrarily dependent p-values
if the p-values are presented to the analyst in a random order.
\end{remark}

We present an extension of the converse direction of Theorem 2.2, which is valid under exchangeability of the vector of p-values. In particular, to preserve exchangeability and use the result stated in Theorem \ref{th:EMI}, it is necessary that the same calibrator $f$ is used for all p-values.
\begin{theorem}\label{th:exch}
    Let $f$ be a calibrator, and $\mathbf{P} =(P_1, \dots, P_K) \in \mathcal{U}^K$ be exchangeable. 
    For each $\alpha \in (0,1)$, we have 
    \[
    \pr\left(\exists k \leq K: \frac{1}{k} \sum_{i=1}^k f\left(\frac{P_i}{\alpha} \right) \geq 1 \right) \leq \alpha.
    \]
\end{theorem}

We first give a formal definition of ex-p-merging function, which is a function that yields a p-value if its argument is a vector of exchangeable p-values.

\begin{definition}
An \emph{ex-p-merging function} is an increasing Borel function $F:[0,1]^{K}\to [0,1]$
such that 
$\pr(F(\mathbf P) \le \alpha)\le \alpha$ for all $\alpha \in (0,1)$ and $\mathbf P \in \mathcal U^K$ that is exchangeable. 
It is \emph{homogeneous} if $F(\gamma \mathbf p)=\gamma  F(\mathbf p)$ for all $ \gamma  \in (0,1]$ and $\mathbf p \in [0,1]^{K}$.
An ex-p-merging function $F$ is \emph{admissible} if for any ex-p-merging function $G$, $G\le F$ implies $G=F$. 
\end{definition}

We now use the result given in Theorem \ref{th:exch} to derive better p-merging functions by exploiting the duality between rejection regions and p-merging functions. In particular, to derive an ex-p-merging the following steps are involved: Initially, the rejection regions at level $\alpha$ based on a given calibrator are determined. As explained in Section \ref{sec:notation}, the ex-p-merging is established by choosing the smallest $\alpha$ for which the p-values $\mathbf{p}$ falls within the rejection region $R_\alpha$. In the first step, the result stated in Theorem \ref{th:exch} helps to derive functions that dominate their counterpart valid under arbitrary dependence. To elaborate, 
starting from a calibrator $f$ and $\alpha \in (0,1)$, we define the {exchangeable rejection region}
\[
R_\alpha = \left\{\mathbf{p} \in [0,1]^K : \frac{1}{k} \sum_{i=1}^k f\left(\frac{p_i}{\alpha} \right) \geq 1 \mathrm{~for~some~}k\leq K\right\}.
\]
Using $R_\alpha$, we can define the function $F:[0,1]^K \to [0,1]$ by
\begin{equation}\label{eq:f_exch}
\begin{split}
F(\mathbf{p}) &= \inf\{\alpha \in (0,1): \mathbf{p} \in R_\alpha \}\\
&= \inf\left\{\alpha \in (0,1): \exists k \leq K \mathrm{~s.t.~}\frac{1}{k} \sum_{i=1}^k f\left(\frac{p_i}{\alpha} \right) \geq 1 \right\}\\
&= \inf\left\{\alpha \in (0,1): \bigvee_{k=1}^K \frac{1}{k} \sum_{i=1}^k f\left( \frac{p_i}{\alpha} \right) \geq 1 \right\},
\end{split}
\end{equation}
where, and throughout, $ \bigvee_{k=1}^K \frac{1}{k} \sum_{i=1}^k f( \frac{p_i}{\alpha} )  $
should be understood as 
$ \bigvee_{k=1}^K (\frac{1}{k} \sum_{i=1}^k f( \frac{p_i}{\alpha} )) $. 
Note that 
\eqref{eq:f_exch} is always smaller or equal than the p-merging function given by \eqref{eq:calibrator}
\[
F'(\mathbf{p}) = \inf\left\{\alpha \in (0,1): \frac{1}{K} \sum_{k=1}^K f\left(\frac{p_k}{\alpha}\right) \geq 1 \right\},
\]
which is valid for p-values with an arbitrary dependence. This is particularly important since all admissible homogeneous and symmetric p-merging functions have the form $F'$ for some admissible calibrator (a symmetric version of Theorem \ref{th:vovkwang}; see \cite{vovk2022}). This implies that the function defined in \eqref{eq:f_exch} dominates the function $F'$. In the following theorem we prove that \eqref{eq:f_exch} is a homogeneous ex-p-merging function.

\begin{theorem}\label{th:calibr_to_F_exc}
    If $f$ is a calibrator and $\mathbf{P}\in \mathcal{U}^K$ is an exchangeable sequence, then $F$ in~\eqref{eq:f_exch} is a homogeneous ex-p-merging function.
\end{theorem}
It is clear from~\eqref{eq:f_exch} that the function depends on the order of the values in $\mathbf{p}$, and hence $F(\mathbf{p})$ is not a symmetric function.  
This is not a coincidence: in the next result, we show that any symmetric  ex-p-merging function is actually valid under arbitrary dependence, and hence they cannot improve over the admissible p-merging functions studied by \cite{vovk2022}. In particular, this implies that under exchangeability the multiplicative factor 2 for the arithmetic average cannot be improved, as earlier noted by \cite{choi2023}, but it also extends their result to every other symmetric merging function.

\begin{proposition}\label{prop:ex-p-m}
A symmetric ex-p-merging function is necessarily a p-merging function. Hence, for an ex-p-merging function to strictly dominate an admissible p-merging function, it cannot be symmetric. 
\end{proposition}

Clearly, Proposition~\ref{prop:ex-p-m} implies that under symmetry,
a function is a p-merging function if and only if it is an \emph{ex-p-merging function}. 
Hence, to take advantage of the exchangeability of the p-values (over arbitrary dependence), one necessarily deviates from symmetric ways of merging p-values, as done in Theorem~\ref{th:calibr_to_F_exc}.
More importantly, the proof of Proposition~\ref{prop:ex-p-m} illustrates the idea (mentioned in Remark~\ref{rem:31}) that for arbitrarily dependent p-values, we can first randomly permute them and then apply an ex-p-merging function (not necessarily symmetric) such as the one in Theorem~\ref{th:calibr_to_F_exc}, to obtain a p-value.

The next result gives a simple condition on the calibrator $f$ that guarantees that the probability of rejection using $F$ in \eqref{eq:f_exch} is sharp for some $\mathbf P$. 

\begin{proposition}
\label{prop:sharp}
Suppose that $f$ is a convex admissible calibrator with $f(0+)\le K$ and $f(1)=0$, and $F$ is in \eqref{eq:f_exch}.  
      For $\alpha \in (0,1)$, there exists an exchangeable $\mathbf P\in \mathcal U^K$ such that $\pr  (F(\mathbf{P}) \leq \alpha   )
      =\alpha. $  
\end{proposition}

\begin{remark}
Proposition~\ref{prop:sharp} is not sufficient to justify admissibility of $F$ in \eqref{eq:f_exch}. 
In general, admissibility of \emph{ex-p-merging functions} remains unclear. 
For instance, take $F$ in \eqref{eq:f_exch} with $f(p)=(2-2p)_+$, corresponding to the arithmetic average, as in  Section~\ref{sec:avg} below.
If $K=2$,  
then $F$ is not admissible as it is strictly dominated by the Bonferroni p-merging function given by $F_{\mathrm{Bonf}}(p_1,\dots,p_K)= K \min\{p_1,\dots,p_K\}$. 
For $K\ge 3$, $F$  and $F_{\mathrm{Bonf}}$ are not comparable. 
\end{remark}

\subsection{Homogeneity of p-merging functions}
 
{As seen from Theorem \ref{th:calibr_to_F_exc}, the class of ex-p-merging functions we obtained in \eqref{eq:f_exch} are homogeneous. Indeed, all explicit p-merging functions in the literature are homogeneous; see \cite{vovk2022} for many examples.
In the next result, we show a very special feature of  homogeneous p-merging functions, justifying their relevance in applications.  

\begin{theorem}
\label{th:r1-1}
Let $F$ be a p-merging function. For any $\alpha \in (0,1)$, there exists a homogeneous p-merging function $G$
such that  $R_\alpha(F) \subseteq R_\alpha(G)$.
\end{theorem}

As a consequence of Theorem \ref{th:r1-1}, if the level $\alpha$ is determined before choosing the p-merging function, then it suffices to consider homogeneous ones, since their rejection sets are at least as larger as  those of other p-merging functions.  
Since all admissible homogeneous p-merging functions have the form in Theorem \ref{th:vovkwang}, the class of our ex-p-merging functions  in \eqref{eq:f_exch}, sharing a similar form to \eqref{eq:calibrator}, is quite broad. 
Note that Theorem \ref{th:r1-1} does not imply that there exists a homogeneous p-merging function $G$ dominating $F$ in general, because the construction of $G$ depends on the given $\alpha.$
}

\subsection{Sequentially combining a stream of exchangeable p-values}
In Subsection~\ref{sec:exch_gen} we have seen that, if our starting vector of p-values $\mathbf{P}=(P_1,\dots,P_K)$ is exchangeable, then it is possible to derive new combination rules exploiting their exchangeability. The technique for developing these new rules involves converting the initial p-values into e-values through an $\alpha$-dependent calibrator. Following this transformation, the exchangeable Markov inequality (Theorem~\ref{th:EMI}) plays a crucial role in formulating these rules. In particular, notice that the inequality in Theorem~\ref{th:EMI} is uniformly valid; therefore, the exchangeable sequence of p-values need not be limited to a set with cardinality $K$, but it is possible to continue to add p-values and stop when the procedure is stable.

To provide a concrete example, in the case where p-values are obtained by applying a sample-splitting procedure to the same dataset, a researcher can obtain one p-value at a time simply by performing a new data split. The researcher would then aim to combine these p-values sequentially and potentially stop when the procedure appears to stabilize.

We first define the following lemma:
\begin{lemma}\label{lemma:continuous}
    Let $\mathbf{p}=(\mathbf{p}_1,\mathbf{p}_2) \in [0,1]^K$ be a vector with $\mathbf{p}_1   \in [0,1]^{K_1}, \mathbf{p}_2  \in [0,1]^{K_2}$ and $K=K_1+K_2$.
    In addition, let $F$ be the function defined in \eqref{eq:f_exch}. Then
    \[
    F(\mathbf p) = F(\mathbf p_1, \mathbf p_2) \leq F(\mathbf p_1).
    \]
\end{lemma}

The above result  implies that the function $F$ in \eqref{eq:f_exch} is non-increasing as more parameters (i.e., p-values) are added. In addition, the following holds.
\begin{theorem}\label{th:sequential}
    Let $P_1, P_2, ... \in \mathcal{U}^\infty$ be an infinite exchangeable sequence and let $F$ be the function defined in \eqref{eq:f_exch}. Then,
    \[
    \pr(\exists k \geq 1: F(\mathbf{P}_k) \leq \alpha) \leq \alpha,
    \]
    where $\mathbf{P}_k = (P_1, \dots, P_k)$. 
\end{theorem}

The result allows the analyst to either continue collecting new (exchangeable) p-values or to cease based on the outcome. In particular, in the example described at the beginning of the section, a researcher can stop the procedure when the result seems to stabilize. However, there are some issues to consider when applying such a procedure. In particular, some calibrators $f$ depend on the number $K$ of p-values. For example, this is the case for the Hommel combination (Section~\ref{sec:hommel}) and for the harmonic mean (Section~\ref{sec:harm_avg}). 
As a final caveat, exchangeability becomes  more stringent when the number $K$ grows; for instance, in general, for a given $K$-dimensional exchangeable vector $(P_1,\dots,P_K)$, there may not exist $P_{K+1}$ such that  $(P_1,\dots,P_{K+1})$ is exchangeable. Luckily, in many practical situations, the procedure that produced the $K$ exchangeable p-values in the first place could also produce more of them; for example, this happens when the p-value was produced by sample splitting.

\subsection{Randomized p-merging functions}\label{sec:rand_gen}
In this subsection, we start with a collection of arbitrarily dependent p-values and we will show how it is possible to enhance existing merging rules using a simple randomization trick. In this case, we denote
\[
\mathcal U^K\otimes \mathcal U=\{(\mathbf P,U)\in  \mathcal U^K\times \mathcal U: \mbox{$U$ and $\mathbf P$ are independent} \}, 
\]
and we state a randomized version of the converse direction of Theorem~\ref{th:vovkwang}, by changing the constant $1$ in \eqref{eq:calibrator} to a uniform random variable $U$.

\begin{theorem}\label{th:1}
Let $f_1,\dots,f_K$ be calibrators and
 $ (P_1,\dots,P_K,U) \in \mathcal U^K\otimes \mathcal U$. For each $\alpha \in (0,1)$ and $(\lambda_1,\dots,\lambda_K)\in\Delta_K$, we have  \begin{equation*}
\pr \left(
      \sum_{k=1}^K \lambda_k f_k\left(\frac{P_k}{\alpha}\right) \ge U
    \right)  \le \alpha.
  \end{equation*} 
  If $ f_1,\dots,f_K$ are  admissible calibrators 
  and $\pr(\sum_{k=1}^K \lambda_k f_k(P_k/\alpha) \le 1)=1$,
  then  equality holds 
    \begin{equation}\label{eq:randomized-2}
\pr   
    \left(
      \sum_{k=1}^K \lambda_k f_k\left(\frac{P_k}{\beta}\right) \ge U
    \right)  = \beta \mbox{~~for all $\beta \in (0,\alpha]$}.
  \end{equation} 
\end{theorem}

The result in Theorem~\ref{th:1} is a direct consequence of the uniformly randomized Markov inequality (UMI) introduced by~\cite{ramdas2023}; see Theorem~\ref{th:UMI}.

\begin{definition}
A \emph{randomized p-merging function} is an increasing Borel function $F:[0,1]^{K+1}\to [0,1]$
such that 
$\pr(F(\mathbf P,U) \le \alpha)\le \alpha$ for all $\alpha \in (0,1)$ and $ (\mathbf P,U) \in \mathcal U^K \otimes \mathcal U$.
It is \emph{homogeneous} if $F(\gamma \mathbf p,u)=\gamma  F(\mathbf p,u)$ for all $ \gamma  \in (0,1]$ and $(\mathbf p,u)\in [0,1]^{K+1}$.
A  randomized p-merging function $F$ is admissible if for any randomized p-merging function $G$, $G\le F$ implies $G=F$. 
\end{definition}

Let $f_1,\dots,f_K$ be calibrators  and $(\lambda_1,\dots,\lambda_K)\in\Delta_K$.
For $\alpha \in (0,1)$, define the randomized rejection region by
$$ 
    R_\alpha 
    = 
    \left\{
      (\mathbf p, u)\in [0,1]^{K+1}:
      \sum_{k=1}^K \lambda_k f_k\left(\frac{p_k}{\alpha}\right) \ge u
    \right\} 
$$
where we set $f_k(p_k/u)=0$ if $u=0$. 
Using $R_{\alpha}$, we can define the function $F:[0,1]^{K+1}\to [0,1]$ by  
\begin{align}
\label{eq:def-f}
F(\mathbf p, u) &= \inf\{\alpha \in (0,1):   (\mathbf p, u)\in R_\alpha \} \notag
\\&=  \inf\left\{\alpha \in (0,1):   \sum_{k=1}^K \lambda_k f_k\left(\frac{p_k}{\alpha}\right) \ge u \right\},
\end{align}
with the convention $0\times \infty=\infty$ (this guarantees $F(\mathbf p,u)=0$ when any component of $(\mathbf p,u)$ is $0$).
\begin{theorem}\label{th:2}
If $f_1,\dots,f_K$ are calibrators and $(\lambda_1,\dots,\lambda_K)\in\Delta_K$, then  $F$   in \eqref{eq:def-f} is a homogeneous randomized p-merging function.
Moreover,  $F$  is lower semicontinuous.
\end{theorem}
 
In case of symmetric p-merging functions (i.e., $F(\mathbf p, u)=F(\mathbf q, u)$ for any permutation $\mathbf q$ of $\mathbf p$), we have the following corollary, which directly follows from Theorem~\ref{th:2}.
\begin{corollary}\label{co:1}
For any calibrator $f$,
\begin{equation*}
F(\mathbf p, u) =  \inf\left\{\alpha \in (0,1):   \frac 1 K\sum_{k=1}^K  f \left(\frac{p_k}{\alpha}\right) \ge  u \right\},
 \end{equation*}
is a homogeneous, symmetric randomized p-merging function.
\end{corollary} 

A simple observation is that replacing $f(p)$ with $f(p)\wedge K$ does not change the function $F$. This observation allows us to only consider calibrators that are bounded above by $K$, which can improve some existing p-merging functions. This is similar to what was done in \cite{vovk2022} in the context of deterministic p-merging functions.

\begin{remark}[P-hacking via repeated derandomization]
    Randomized methods may not cause a lack of reproducibility if they are part of a standard automated data analysis pipeline without a human in the loop. However, if a human is actively involved in the analysis, then there is a risk of ``p-hacking", where a (malevolent) researcher re-runs the randomized method many times in order to obtain a desirable result according to their own utility. In our setting, the issue translates into sampling different $U$ and picking the smallest one (which would be closer to $0$ the more times the procedure is re-run). Clearly, this procedure is not valid and can be particularly problematic in the context of confirmatory analysis, where the end product is a binary decision.
\end{remark}

\begin{remark}[Internal randomization]
The use of \emph{internal} (as opposed to \emph{external}) randomization can be particularly useful to prevent and mitigate the risk of p-hacking mentioned above. We mention two such strategies below.
First, note that $F $ in \eqref{eq:def-f} is increasing in each of its arguments. 
If one has prior information that one of the p-values, say $P_1$, is independent of the rest (but the rest can be arbitrarily dependent), then one can use $P_1$ for randomization and apply $F$ (with one less input dimension for $\mathbf p$) to $(\mathbf p,u)=(P_2,\dots,P_K,U)$ with $U=P_1$ to obtain a p-value that does not depend on external randomization.   
Monotonicity of $ u\mapsto F(\mathbf p,u)$ guarantees two things. 
First,  a p-variable may be stochastically larger than a standard uniform one, so  increasing monotonicity is needed for validity. Second, if $P_1$ is indeed very small, i.e., it carries signal against the null, then the combined p-value will benefit from this signal.  
This form of internal randomization has been discussed in \citet[Section B.1]{wang2024testing}. 

An alternative method of internal randomization through data (instead of p-values) is discussed by \citet[Section 10.6]{ramdas2023}. To understand their proposal, assume that each p-value is calculated using a function of only the order statistics of iid data; for example, it could be based only on sums, like the t-statistic. In this case, the rank of the first data point (amongst all the data points) is a discrete uniform random variable, and can be used in place of $U$. This way, if the dataset is itself public (posted by a previous research paper, for example), the ordering itself is not in the hands of the researcher analyzing that data, reducing the risk of p-hacking. We refer interested readers to \citet[Section 10]{ramdas2023} or \citet{lei2024inference} for further discussions.
\end{remark}

It is feasible to combine the results presented in Subections~\ref{sec:exch_gen} and~\ref{sec:rand_gen} through the formulation of novel p-merging functions that exploit both the properties of exchangeability and randomization. These results are presented in Appendix \ref{sec:exch_rand_pmerg} and are based on the exchangeable and uniformly randomized Markov inequality presented in Theorem~\ref{th:EUMI}.

\subsection{Instantiating the above ideas}

The ideas above were admittedly somewhat abstract, but provide us with the general tools to improve specific combination rules. The following sections do this for several rules, one by one. To elaborate, one of the most commonly employed p-merging functions is the Bonferroni method:
\[
F_{\mathrm{Bonf}}(\mathbf{p}) = K p_{(1)},
\]
where $p_{(1)}$ is the minimum of observed p-values. \cite{ruger1978} extended the aforementioned rule in a more general sense. In particular, it is possible to prove that
\begin{equation} \label{eq:ruger}
    F_\mathrm{R}(\mathbf{p}):=\frac{K}{k}p_{(k)}, \quad k\in\{1, \dots, K\},
\end{equation}
is a p-value, where $p_{(k)}$ represents the $k$-th smallest p-value among $(p_1, \dots, p_K)$. In other words, the $\lambda$-quantile $p_{(\lceil \lambda K \rceil)}$ is a p-value if multiplied by the factor $1/\lambda$. In particular, the robust and widely used combination rule, twice the median, is part of this class. The next section improves on this combination rule.

 \cite{vovk2020} introduced the class of p-merging functions based on the generalized mean, also called \emph{M-family}. This general class takes the form
\begin{equation}\label{eq:gen_mean}
   a_{r,K}\left(\frac{p_1^r + \dots + p_K^r}{K}\right)^{1/r},
\end{equation}
where $r \in \mathbb{R} \setminus \{0\}$ and $a_{r,K}$ is the smallest constant making \eqref{eq:gen_mean} a p-merging function. This class encompasses numerous well-known cases, each distinguished by different values of the parameter $r$. In particular, if $r=1$ then~\eqref{eq:gen_mean} reduces to the simple average introduced in \cite{ruschendorf1982} and the value $a_{1,K}=2$. Another important case is the harmonic mean obtained with $r=-1$. Among this class, the harmonic mean combination rule should be used when substantial dependence among the p-values is suspected. If the dependence is really strong, the arithmetic mean might be a safer option.
In the following sections, we demonstrate that if p-values exhibit exchangeability or if randomization is allowed, then it becomes feasible to enhance most of these combination rules. 

Before continuing, we provide a few simple examples to highlight the benefits of the proposed findings and demonstrate how our methods can enhance the existing approaches. In the examples we will use some rules proposed in Table \ref{tab:1}, foreshadowing many results to come.

\begin{example}
    
Suppose that the vector  $\mathbf{P}$ of $3$  p-values is generated as follows. 
With $0.9$ probability, 
$\mathbf{P}=(P_1,P_2,P_3)$
where $P_1,P_2,P_3$ are independent, and 
with $0.1$ probability $\mathbf{P}=(P_4,P_4,P_4) $. 
We assume $P_i \in \mathcal{U}$ under the null, while under the alternative each   $P_i$ is distributed as  $\mathrm{Beta}(0.2,1)$.
The Beta distribution $(a,1)$ with small $a>0$ is a typical model for p-values under the alternative hypothesis; see 
\cite{sellke2001calibration}. 
Clearly, the p-values are exchangeable under the null. Suppose that one want to use the rule $(3/2)p_{(2)}$ to combine the p-values. Then we can check numerically that, fixing the threshold $\alpha$ to $0.05$, the probability of rejection is $0.5101$ under the alternative. If we use the ex-p-merging derived from the median (see Theorem \ref{thm:ruger_exch} below) the   probability of rejection increases to $0.6207$.
\end{example}

\begin{example}
Suppose that we want to test an hypothesis and we have that under the null $P_1 \in \mathcal{U}$ while under the alternative $P_1 \sim \mathrm{Beta(0.2,1)}$. The vector of p-values is generated as follows $\mathbf{P}:=(P_1,P_2)=(P_1,1-P_1)$, where $P_2$ is an exact p-value and $\mathbf{P}$ is exchangeable under the null. In this scenario, the commonly applied twice the mean, even though controls the type I error, has no power under the alternative hypothesis because the result is always equal to $1$. On the other hand, using the exchangeable rule derived from the arithmetic average (see Theorem \ref{thm:avg_exch}),  the probability of rejection is $\pr(P_1 \leq \alpha/2) \approx 0.48$ under the alternative for $\alpha = 0.05$.
\end{example}

\subsection{Combining asymptotic p-values}

Before proceeding with the remainder of the paper and introducing new merging rules based on the results presented in the preceding sections, we want to examine the scenario wherein the p-values are asymptotically valid. Many of the results obtained in the literature rely on uniform or super-uniform p-values (see Section~\ref{sec:notation}); however, in statistical applications, p-values are often asymptotic, and they are not necessarily valid p-values in finite sample. See, for example, \cite{severini2000likelihood} for an introduction to p-values obtained using likelihood-based methods.

All methods in our paper work also for asymptotic p-values, i.e., those that converge in distribution to p-values. 
Suppose that $\mathbf (\mathbf P_n)_{n\in \mathbb N}$
is a sequence of nonnegative random vectors that converges to a vector $\mathbf P$ of p-values in distribution.
With an upper semicontinuous calibrator $f$ (recall that all admissible calibrators are upper semicontinuous),
for each $\alpha\in (0,1)$ and $u\in (0,1)$,
the rejection sets  $R_\alpha$ given by 
$$R_\alpha=\left\{ \mathbf{p} \in [0,\infty)^K: \bigvee_{k \le K} \frac{1}{k} \sum_{i=1}^k f\left(\frac{p_i}{\alpha} \right) \geq 1 \right\}
$$
or 
$$R_{\alpha}=\left\{\mathbf{p} \in [0,\infty)^K:   \frac 1 K\sum_{k=1}^K  f \left(\frac{p_k}{\alpha}\right) \ge  u \right\}$$
are closed. 
As a consequence, by the Portmanteau Theorem, 
$$\limsup_{n\to \infty} \mathbb P(\mathbf P_n \in R_{\alpha} )\le \mathbb P (\mathbf P \in R_{\alpha} ) .$$
Therefore, any methods  in our paper that produce a p-value for the vector $\mathbf P$ of p-values (exchangeable or arbitrarily dependent)
produce an asymptotic p-value for  any $(\mathbf P_n)_{n\in \mathbb N}$ that converges to $\mathbf P$ in distribution.

\section{Improving R\"uger's combination rule}\label{sec:ruger}

\cite{vovk2022} showed that the p-merging function defined in~\eqref{eq:ruger}, with a trivial modification (i.e., return $0$ if $p_{(1)}=0$; see Theorem 7.3 of \cite{vovk2022}) is admissible for $k\ne K$, and it is admissible among symmetric p-merging functions when $k=K$. In particular, the corresponding calibrator that induces~\eqref{eq:ruger} is $$f(p) = \frac{K}{k} \ind\{p \in (0, k/K]\} + \infty \ind\{p=0 \},$$ and this implies that we can exploit directly the duality between rejection regions and p-merging functions.

\subsection{An exchangeable R\"uger combination rule}
If the exchangeability condition is satisfied, then we can obtain something better than the combination rule in~\eqref{eq:ruger}.  First, it is clear that $f(p) = \frac{K}{k} \ind\{p \in (0, k/K]\} + \infty \ind\{p=0 \}$ is an admissible calibrator. We now define the function
\begin{equation*}
    F_{\mathrm{ER}}(\mathbf{p}) = \inf\left\{\alpha \in (0,1): \bigvee_{\ell \le K} \frac{1}{\ell} \sum_{i=1}^\ell \frac{K}{k} \ind\left\{\frac{p_i}{\alpha} \leq \frac{k}{K} \right\} \geq 1 \right\}.
\end{equation*}
Below, for fixed $k\in \{1,\dots,K\}$, we let $p^\ell_{(\lambda_\ell)}$ denote the $\lceil \ell \frac{k}{K} \rceil$-th ordered value obtained using the first $\ell$ values of $\mathbf{p}$. Essentially, $p^\ell_{(\lambda_\ell)}$ is the upper quantile of order $k/K$ obtained using the first $\ell$ p-values.

\begin{theorem}\label{thm:ruger_exch}
For any fixed $k\in\{1, \dots, K\}$ the function $F_{\mathrm{ER}}$ satisfies
\[
F_{\mathrm{ER}}(\mathbf{p}) = \left(\frac{K}{k} \bigwedge_{\ell=1}^K p^\ell_{(\lambda_\ell)}\right) \ind\{p_{(1)} > 0\} \mbox{~for  $ \mathbf{p} \in [0,1]^K$, }
\]
where $\lambda_\ell := \lceil \ell \frac{k}{K} \rceil$, and it is an ex-p-merging function that strictly dominates $F_{\mathrm{R}}$ in~\eqref{eq:ruger}.
\end{theorem}

It may be useful to note that the Bonferroni rule is not improved using this method. Indeed, fixing $k=1$, we find that $p_{(\lambda_\ell)}^\ell$ reduces to the minimum of the first $\ell$ p-values subsequently taking the minimum of the obtained sequences coincides with the overall minimum. In addition, R\"uger can be sharp, for some exchangeable $\mathbf{P}$, i.e., satisfying $F_{\rm R}(\mathbf P)\in \mathcal U$, and so is our proposal.

\subsection{A randomized R\"uger combination rule}
In this part, we prove that if randomization is allowed, it becomes feasible to enhance the combination introduced by~\cite{ruger1978}, even if the sequence of p-values presents an arbitrary dependence and the obtained result has nice properties in terms of interpretability. As before, we define the merging function
\begin{equation*}
F_{\mathrm{UR}}(\mathbf{p},u) = \inf\left\{\alpha \in (0,1): \frac{1}{K} \sum_{i=1}^K \frac{K}{k} \ind\left\{\frac{p_i}{\alpha} \leq \frac{k}{K} \right\} \geq u \right\}.
\end{equation*}

\begin{theorem}\label{th:ruger_rand}
    For any fixed $k \in \{1, \dots, K\}$, the function $F_{\mathrm{UR}}$ satisfies
    \[
    F_{\mathrm{UR}}(\mathbf{p}, u) = \frac{K}{k} p_{(\lceil uk \rceil)}\, \ind\{p_{(1)} > 0\},
    \]
    and it is a randomized p-merging function that strictly dominates $F_{\mathrm{R}}$  in~\eqref{eq:ruger}.
\end{theorem}

The above theorem implies in particular that the R\"uger combination rule is inadmissible if external randomization is allowed, despite it being admissible if randomization is not allowed~\citep{vovk2022}. The fact that $p_{(\lceil {UK} \rceil)}$ is a p-value is particularly interesting. It has a simple interpretation: sort the p-values and pick the one at a uniformly random index. In addition, for the latter combination rule \eqref{eq:randomized-2} holds for all $\beta \in (0,1]$.

It is worth noting that when $k=1$ the R\"uger combination rule reduces to the Bonferroni method; however, the introduction of a randomization does not yield any practical benefit since $\lceil {U} \rceil = 1$.   

\section{Improving Hommel's combination rule}\label{sec:hommel}
The method proposed in Section~\ref{sec:ruger} requires us to choose the value of $k$ in advance; a solution that solves the problem is proposed by~\cite{hommel1983}. Hommel's combination rule is given by
\begin{equation}\label{eq:hommel_cr}
    F'_{\mathrm{Hom}}(\mathbf{p}) := h_K \bigwedge_{k=1}^K F_{\mathrm{R}}(\mathbf{p};k) = \left(\sum_{k=1}^K \frac{1}{k}\right)  \bigwedge_{k=1}^K \frac{K}{k} p_{(k)}, 
\end{equation} 
with $h_K= \sum_{k=1}^K \frac{1}{k}$. This function allows selecting the minimum derived from the combinations based on ordered statistics with a multiplicative cost of $h_K \approx \log K$.

It is possible to prove that the Hommel combination rule is not admissible and is dominated by the \emph{grid harmonic merging function} introduced in~\cite{vovk2022}. For completeness, we state here a useful lemma.
\begin{lemma}\label{lemma:calib_hommel}
    Let $f$ be a function defined by
    \begin{equation}\label{eq:grid_harm_calib}
        f(p) = \frac{K \ind\{h_K p \leq 1\}}{\lceil K h_K p\rceil}.
    \end{equation}
    Then $f$ is an admissible calibrator. Moreover, the p-merging function induced by $f$ is
    \begin{equation*}
            F_\mathrm{Hom}(\mathbf{p}):= \inf\left\{\alpha \in (0,1): \sum_{k=1}^K \frac{\ind\{h_K p_k/\alpha \leq 1\}}{\lceil K h_K p_k/\alpha\rceil} \geq 1 \right\},
    \end{equation*}
    is valid and it dominates the Hommel combination rule.
\end{lemma}
The calibrator in~\eqref{eq:grid_harm_calib} coincides, with a slight adjustment, with the BY calibrator introduced in~\cite{xu2022post}.
An interesting fact is that the function $F_{\mathrm{Hom}}$ is always admissible in the class of symmetric p-merging functions; while it is admissible in the class of p-merging function (not necessarily symmetric) if $K$ is not a prime number \citep[Theorem 7.1]{vovk2022}. The Hommel function allows for the selection of the minimum among the $K$ possible different quantiles of $\mathbf{p}$.

In reality, one can choose to select only certain quantiles among $K$ (e.g., one can select the minimum between $K$ times the minimum, 2 times the median and the maximum), hoping to pay a price less than $\log K$. We treat this problem in Appendix \ref{sec:app_genHom}, where we introduce a generalization of the Hommel combination rule. 

\subsection{An exchangeable Hommel's combination rule}
Starting from the results in the previous sections, we can introduce a merging function, which improves Hommel's combination if the input p-values are exchangeable. In particular, we define the function
\begin{equation}\label{eq:feh*}
    F_{\mathrm{EHom}}(\mathbf{p}) = \inf\left\{\alpha \in (0,1): \bigvee_{\ell \le K} \frac{1}{\ell}\sum_{i=1}^\ell \frac{K \ind\{h_K p_i/\alpha \leq 1\}}{\lceil K h_K p_i/\alpha \rceil} \geq 1 \right\}.
\end{equation}
\begin{theorem}\label{thm:hom_exch}
    The function $F_{\mathrm{EHom}}$ is an ex-p-merging function and it strictly dominates the function $F_{\mathrm{Hom}}$.
\end{theorem}

The computation of a closed form for~\eqref{eq:feh*} is not straightforward, a possible solution to calculate the value of $F_{\mathrm{EHom}}$ is by using Algorithm~\ref{alg:alpha_search} defined in Appendix \ref{subsec:algo}. 

\subsection{A randomized Hommel's combination rule}
The \emph{randomized} version of the function $F_{\mathrm{Hom}}$ has been proposed in~\citet[Appendix E]{xu2023} and takes the following form:
\begin{equation*}
    F_{\mathrm{UHom}}(\mathbf{p}, u) := \inf\left\{\alpha \in (0,1): \sum_{k=1}^K \frac{\ind\{h_K p_k/\alpha \leq 1\}}{\lceil K h_K p_k/\alpha\rceil} \geq u \right\}.
\end{equation*}
For completeness, we report here the following theorem.
\begin{theorem}\label{thm:hom_rand}
    The function $F_{\mathrm{UHom}}$ is a randomized p-merging function and it strictly dominates the function $F_{\mathrm{Hom}}$.
\end{theorem}
The value of the function $F_{\mathrm{UHom}}$ can be computed using Algorithm 1 in~\cite{vovk2022}, substituting the value of $1$ by $u$.

\section{Improving the ``twice the average'' combination rule}\label{sec:avg}
We now study the case of $r=1$ in \eqref{eq:gen_mean}, which corresponds to the arithmetic mean. The general case, which allows $r \in \mathbb{R} \setminus \{0\}$, is considered in Appendix \ref{sec:gen_mean}.
In the following, let $A(\mathbf{p})$ denote the arithmetic average of any vector $\mathbf{p}$, let $\mathbf{p}_{m} := (p_1,\dots,p_m)$ denote the vector containing the first $m$ p-values, and let $\mathbf{p}_{(m)}$ denote the vector containing the smallest $m$ elements of $\mathbf{p}$: $\mathbf{p}_{(m)} = (p_{(1)}, \dots, p_{(m)})$ such that $p_{(1)}\leq \dots \leq p_{(m)}$. In addition, we denote by $\mathbf{p}_{(m)}^\ell = (p_{(1)}^\ell, \dots, p_{(m)}^\ell)$, $m \in \{1, \dots, \ell\}$, the vector containing the smallest $m$ elements of $(p_1, \dots, p_\ell)$.
First, let us introduce a lemma that will be instrumental in subsequent discussions. 
\begin{lemma}\label{lemma:calibrator}
    Let $f(p) = (2- 2p)_+\ind\{p \in (0,1]\} + \infty \ind\{p=0\}$. Then, $f$ is an admissible calibrator.
\end{lemma}
In particular, we have that the calibrator defined in Lemma~\ref{lemma:calibrator} is larger or equal than $f'(p) := (2- 2p)$, which is the function inducing the average combination rule  
\begin{equation}\label{eq:merging-avg}
F'_{\mathrm{A}}(\mathbf{p}) := 2 A(\mathbf{p}),
\end{equation}
which is a valid p-merging function.
Note that $f'$ can take negative values (its arguments $p/\alpha$ can be larger than $1$), so it is technically not a calibrator in the sense of our definitions.

\subsection{Exchangeable average combination rule}
We now define  ex-p-merging functions
\begin{align*}
    F_{\mathrm{EA}}(\mathbf{p}) &= \inf\left\{\alpha \in (0,1): \bigvee_{\ell \le K} \frac{1}{\ell} \sum_{i=1}^\ell \left(2 - 2 \frac{p_i}{\alpha} \right)_+ \geq 1 \right\};\\
    F_{\mathrm{EA}}'(\mathbf{p}) &= \inf\left\{\alpha \in (0,1): \bigvee_{\ell \le K} \frac{1}{\ell} \sum_{i=1}^\ell \left(2 - 2 \frac{p_i}{\alpha} \right) \geq 1 \right\}.
\end{align*}

\begin{theorem}\label{thm:avg_exch}
  The dominations among ex-p-merging functions $F_{\mathrm{EA}}\le F'_{\mathrm{EA}}\le F_\mathrm{A}'$ are strict. Moreover, 
    \begin{align*}
    F'_{\mathrm{EA}} (\mathbf p)&= 2\left\{ \bigwedge_{m=1}^K A(\mathbf{p}_{m}) \right\} ;\\
    F_{\mathrm{EA}} (\mathbf p)&= \left( \bigwedge_{\ell=1}^K  \bigwedge_{m=1}^\ell \frac{2A(\mathbf{p}_{(m)}^\ell)}{(2-\ell/m)_+}  \right) \, \ind\{p_{(1)} > 0\}.
    \end{align*}
\end{theorem}

Despite being strictly dominated, $F'_{\mathrm{EA}}$ is very interpretable: it is just the minimum (over $m$) of ``twice the average'' of the first $m$ p-values.

\subsection{Randomized average combination rule}\label{sec:rand_avg}
We now derive an improvement for the ``twice the average" rule using a simple randomization trick. In this case, we do not require exchangeability but we allow for an arbitrary dependence among the p-values. We define the randomized p-merging function $F_{\mathrm{UA}}$ by
\begin{equation*}
F_{\mathrm{UA}}(\mathbf p, u) =  \inf\left\{\alpha \in (0,1):   \frac{1}{K}\sum_{k=1}^K   \left( 2-\frac{2 p_k}{\alpha}\right)_+ \ge u \right\}.
\end{equation*}
Clearly, $F_{\mathrm{UA}}(\mathbf p,u) \le F'_{\mathrm{UA}}(\mathbf p,u)$, where $F'_{\mathrm{UA}}$ is defined by
\begin{equation*}
   F'_{\mathrm{UA}}(\mathbf p, u)  
 =  \inf\left\{\alpha \in (0,1):   \frac{1}{K}\sum_{k=1}^K   \left( 2-\frac{2 p_k}{\alpha}\right) \ge u \right\} .
\end{equation*}
In particular, if randomization is not allowed then $u$ is replaced by $1$ and $F'_{\mathrm{UA}}(\mathbf{p}, 1)$ coincides with $F_\mathrm{A}'$ in~\eqref{eq:merging-avg}.
A p-merging function (such as $F_\mathrm{A}'$) can also be seen as a randomized p-merging function, for which the argument $u$ does not affect its value.

\begin{theorem}\label{theorem:FUA}
  The dominations among randomized p-merging functions $F_{\mathrm{UA}}\le F'_{\mathrm{UA}}\le F_\mathrm{A}'$ are strict. Moreover, 
    \begin{align*}
        F'_{\mathrm{UA}} (\mathbf p, u)&= \frac{2A(\mathbf{p})}{2-u};\\
        F_{\mathrm{UA}} (\mathbf p, u)&= \left(\bigwedge_{m=1}^K \frac{2A(\mathbf{p}_{(m)})}{(2-Ku/m)_+}  \right) \ind\{p_{(1)} > 0\}.
    \end{align*}
    
\end{theorem}

A method that can be directly compared with Theorem~\ref{theorem:FUA} is to use 
$ F_{\mathrm{UA}}^*(\mathbf p, u):=A(\mathbf p)/(2-2u)$ proposed by \citet[Section B.2]{wang2024testing}.
This function $  F_{\mathrm{UA}}^*$ is also a randomized p-merging function.
One can see that $ F_{\mathrm{UA}}^*$ 
does not dominate and is not dominated by any of $F_\mathrm{A}$, $F_{\mathrm{UA}}$ and  $ F_{\mathrm{UA}}'$.
Moreover, there is a simple relationship: $ F_{\mathrm{UA}}'(\mathbf p,u) \le  F_{\mathrm{UA}}^*(\mathbf p,u)$ if and only if $u\ge 2/3$ for   every $\mathbf p$ that is not the zero vector.

\section{Improving the harmonic mean combination rule}\label{sec:harm_avg}
The harmonic mean p-value was studied by \cite{wilson2019}. In our context, it corresponds to the merging function  in~\eqref{eq:gen_mean} when $r=-1$.   We first state a lemma on a calibrator that we will use later.
\begin{lemma}\label{lemma:calib_harmonic}
    Define the function
    \[
    f(p) = \min\left\{\frac{1}{T_K p} - \frac{1}{T_K}, K \right\} \ind\{p \in [0,1]\},
    \]
    with $T_K \geq 1$. Then $f$ is a calibrator if $T_K$ satisfies $KT_K+1-e^{T_K} \leq 0$, and in particular, $f$ is a calibrator if $T_K = \log K + \log \log K + 1$.
\end{lemma}
In the following, we fix  $T_K=\log K + \log \log K + 1 $  and denote by $H(\mathbf{p})=K (\sum_{k=1}^K 1/p_k)^{-1}$ the harmonic mean of the vector $\mathbf{p}$ where $K$ is the number of elements contained in $\mathbf{p}$. We begin with a result that is new even under arbitrary dependence.
\begin{proposition}\label{prop:harm_avg}
    $F'_{\mathrm H}(\mathbf{p}) := (T_K+1)H(\mathbf{p})$ is a p-merging function. 
\end{proposition}

The above result differs from the formulation given in~\citet[Proposition 9]{vovk2020}, which states that $e \log K H(\mathbf{p})$ is a p-merging function, thus sharpening their result for $K \geq 4$. Below we will further  improve this result for exchangeable p-values. Moreover,  a correction factor in the order of $\log K$ as $K\to\infty$ (although smaller than $T_K$)  is needed even for independent p-values 
(see Proposition 6 of \cite{chen2024sub}).
The harmonic mean has some advantages over many other rules under certain dependence conditions;  see \citet{gui2023aggregating}.
It performs similarly to  the Hommel combination; see \cite{chen2023trade}.

\subsection{Exchangeable harmonic mean combination rule}
We define  homogeneous ex-p-merging functions as follows
\begin{align*}
    F_{\mathrm{EH}}(\mathbf{p}) &= \inf\left\{\alpha \in (0,1): \bigvee_{\ell \le K} \frac{1}{\ell} \sum_{i=1}^\ell \left(\frac{\alpha}{T_K p_i} - \frac{1}{T_K} \right)_+ \geq 1 \right\};\\ 
    F'_{\mathrm{EH}}(\mathbf{p}) &= \inf\left\{\alpha \in (0,1): \bigvee_{\ell \le K} \frac{1}{\ell} \sum_{i=1}^\ell \left(\frac{\alpha}{T_K p_i} - \frac{1}{T_K} \right) \geq 1 \right\}.  
\end{align*}

\begin{theorem}\label{thm:harm_exc}
    The dominations among ex-p-merging functions $F_{\mathrm{EH}} \le F'_{\mathrm{EH}}\le F'_{\mathrm{H}}$ 
    are strict. Moreover,
    \begin{align*}
    F'_{\mathrm{EH}}(\mathbf{p})&=\bigwedge_{m=1}^K (T_K+1) H(\mathbf{p}_{m});\\
    F_{\mathrm{EH}}(\mathbf{p})&=\bigwedge_{\ell=1}^K \left( \bigwedge_{m=1}^\ell \left(\frac{\ell\,T_K}{m} + 1\right) H(\mathbf{p}_{(m)}^\ell) \right).
    \end{align*}
\end{theorem}

\subsection{Randomized harmonic mean combination rule}
Similarly to Section~\ref{sec:avg}, we derive an improvement for the harmonic mean using a randomization trick in the case of arbitrarily dependent p-values. 
Define
\begin{align*}
    F_{\mathrm{UH}}(\mathbf{p}, u) &= \inf \left\{\alpha \in (0,1): \frac{1}{K} \sum_{k=1}^K \left(\frac{\alpha}{T_K p_k} - \frac{1}{T_K} \right)_+ \geq u\right\};\\
    F'_{\mathrm{UH}}(\mathbf{p}, u) &= \inf \left\{\alpha \in (0,1): \frac{1}{K} \sum_{k=1}^K \left(\frac{\alpha}{T_K p_k} - \frac{1}{T_K} \right) \geq u\right\}.
\end{align*}
\begin{theorem}\label{thm:harm_rand}
The dominations among randomized p-merging functions $F_{\mathrm{UH}}\le F'_{\mathrm{UH}}\le F_\mathrm{H}'$ are strict. Moreover,  
    \begin{align*}
    F'_{\mathrm{UH}}(\mathbf{p},u)&=(T_K u + 1)H(\mathbf{p});\\
    F_{\mathrm{UH}}(\mathbf{p},u)&=\bigwedge_{m=1}^K \left(\frac{uKT_K}{m} + 1 \right) H(\mathbf{p}_{(m)}).
    \end{align*}
\end{theorem}

A non-randomized improvement of $F'_{\mathrm{H}}$ can be achieved fixing $u=1$ in $F_{\mathrm{UH}}$. This coincides with the function $F_\mathrm{H}(\mathbf{p}) = \bigwedge_{m=1}^K ((KT_K)/m + 1) H(\mathbf{p}_{(m)})$.

\section{Improving the geometric mean combination rule}\label{sec:geo_avg}
We now derive some new combination based on the geometric mean, a special case of~\eqref{eq:gen_mean} when $r\to 0$. Let $G(\mathbf{p})=(\prod_{k=1}^K p_k)^{1/K}$ denote the geometric mean of the vector $\mathbf{p}$. The calibrator, in this case, is given by $f(p) = (-\log p)_+$, which is an admissible calibrator. Actually, a slightly improved calibrator is $f(p)=(-(\log p)/T)_+ \wedge K$ for some $T < 1$ satisfying $\int_0^1 f(p) \d p \leq 1$. This condition is verified when $1-e^{-KT} \leq T$, which makes $T$ very close to 1 for $K$ moderately large~\citep[see Section 3.2 in][]{vovk2020}. In the sequel, we denote by
\begin{equation}\label{eq:geom-merg}
F'_{\mathrm{G}}(\mathbf{p}) := eG(\mathbf{p}),    
\end{equation}
which is a valid p-merging function studied by~\cite{vovk2020}.

\subsection{Exchangeable geometric mean combination rule}
Following the same approach as in the preceding sections, we define  

\begin{align*}
    F_{\mathrm{EG}}(\mathbf{p}) &= \inf\left\{\alpha \in (0,1): \bigvee_{\ell \le K} \frac{1}{\ell} \sum_{i=1}^\ell \left(\log \frac{\alpha}{p_i}\right)_+ \geq 1 \right\};\\
    F'_{\mathrm{EG}}(\mathbf{p}) &= \inf\left\{\alpha \in (0,1): \bigvee_{\ell\le K} \frac{1}{\ell} \sum_{i=1}^\ell \log\frac{\alpha}{p_i} \geq 1 \right\}.
\end{align*}
\begin{theorem}\label{thm:geom_exch}
 The dominations among ex-p-merging functions $F_{\mathrm{EG}} \le F'_{\mathrm{EG}}\le F'_{\mathrm{G}}$ 
    are strict. Moreover, 
    \begin{align*}
    F'_{\mathrm{EG}}(\mathbf{p})&= e\left\{ \bigwedge_{m=1}^K G(\mathbf{p}_{m}) \right\} ;\\
    F_{\mathrm{EG}}(\mathbf{p})&= \bigwedge_{\ell=1}^K \left( \bigwedge_{m=1}^\ell e^{\ell/m} G(\mathbf{p}^\ell_{(m)}) \right).
    \end{align*}
\end{theorem}

\subsection{Randomized geometric mean combination rule}
As in the previous sections, we define the randomized p-merging functions as follows:
\begin{align*}
    F_{\mathrm{UG}}(\mathbf{p}, u) &= \inf\left\{\alpha \in (0,1): \frac{1}{K} \sum_{k=1}^K\left(\log\frac{\alpha}{p_k}\right)_+ \geq u \right\};\\
    F'_{\mathrm{UG}}(\mathbf{p}, u) &= \inf\left\{\alpha \in (0,1): \frac{1}{K} \sum_{k=1}^K\log\frac{\alpha}{p_k} \geq u \right\}. 
\end{align*}

\begin{theorem}\label{thm:geom_rand}
The dominations among randomized p-merging functions $F_{\mathrm{UG}}\le F'_{\mathrm{UG}}\le F_\mathrm{G}'$ are strict. Moreover,  
\begin{align*}
    F'_{\mathrm{UG}}(\mathbf p, u)&=
    e^u G(\mathbf{p});\\
    F_{\mathrm{UG}}(\mathbf p, u)&=
    \bigwedge_{m=1}^K \left( e^{u \frac{K}{m}} G(\mathbf{p}_{(m)}) \right).
\end{align*}
    
\end{theorem}

A non-randomized improvement of the combination in~\eqref{eq:geom-merg} can be obtained fixing $u=1$ in $F_{\mathrm{UG}}$. This gives the combination rule $F_\mathrm{G}(\mathbf{p}) = \bigwedge_{m=1}^K e^{K/m} G(\mathbf{p}_{(m)})$.

\section{Simulation study}\label{sec:simul}
In the previous sections, new p-merging functions have been introduced. These new rules are obtained using a randomization trick or they rely on exchangeability of p-values. Specifically, the introduced rules have been shown to dominate their original counterparts by utilizing randomness or exchangeability (or both). In this section, our aim is to investigate their performance using simulated data. 

We consider the example described in~\citet[Section 6]{vovk2020} (a similar example is proposed in~\cite{chen2023trade}), where p-values are generated in the following way:
\begin{equation}\label{eq:sim_pvals}
    X_k = \rho Z + \sqrt{1-\rho^2}Z_k - \mu,\quad P_k = \Phi(X_k),
\end{equation}
where $\Phi(\cdot)$ is the cumulative density function of the standard normal distribution, $Z, Z_1, \dots, Z_K \stackrel{iid}{\sim} \mathcal{N}(0,1)$, and $\mu \geq 0$ and $\rho \in [0,1]$ are constants. It is simple to prove that $P_1, \dots, P_K$ are exchangeable and their marginal distribution does not depend on $\rho$. In addition, if $\rho=0$ then $P_1, \dots, P_K$ are independent while if $\rho = 1$ then $P_1=\dots=P_K$. The value $p_k$, in this case, can be interpreted as the p-value resulting from a one-side z-test of the null hypothesis $\mu=0$ against the alternative $\mu>0$ from the statistic $X_k \sim \mathcal{N}(-\mu, 1)$ with unknown $\mu$. We let the parameter $\mu$ vary in the interval $[0,3]$ (if $\mu =0$ then $H_0$ is true) and fix the upper bound of the type I error to the nominal level $0.05$.

\begin{figure}
    \centering
    \includegraphics[width = 0.9\linewidth]{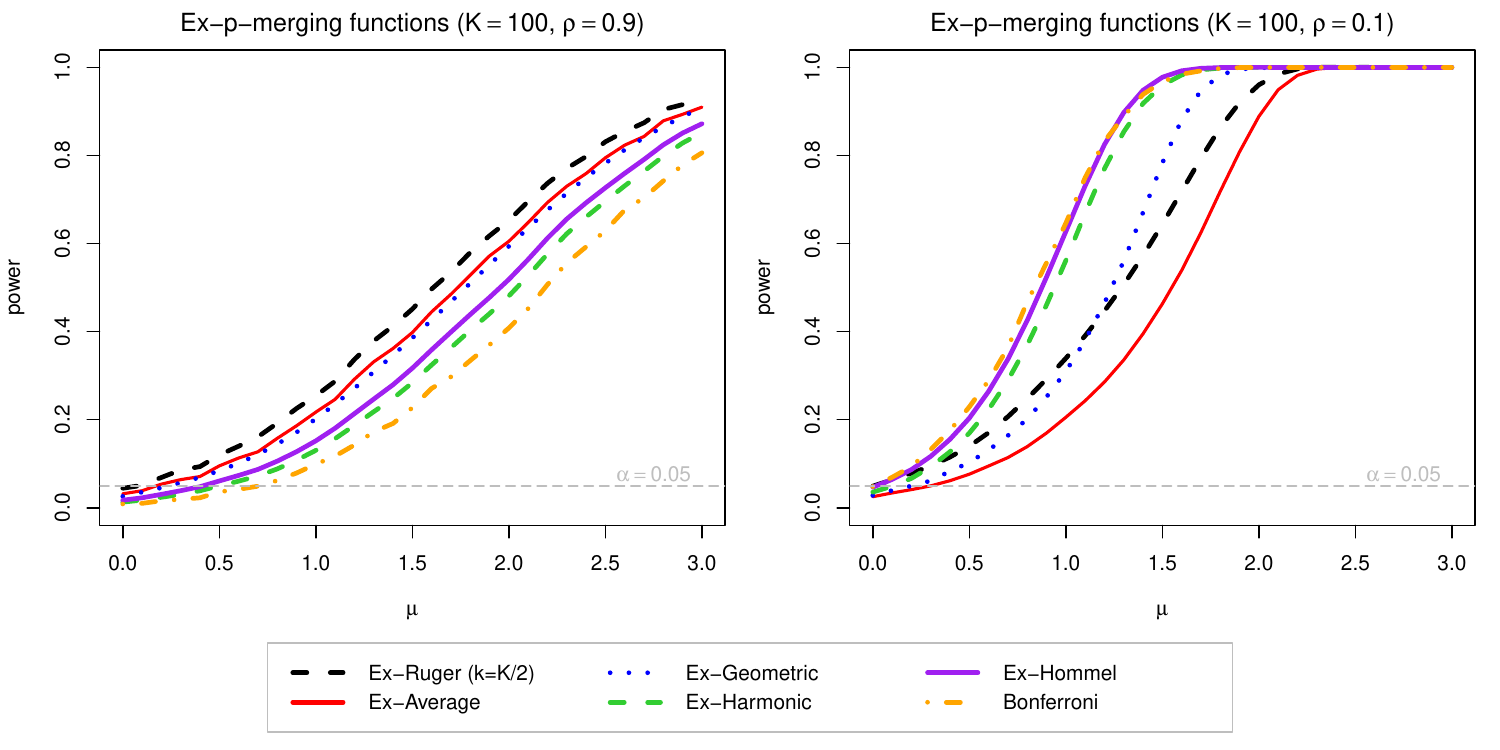}
    \caption{Combination of p-values using different {ex-p-merging functions} under high (left) and low (right) dependence. The performance of the different {ex-p-merging functions} is almost reversed in the two situations.}
    \label{fig:exch_rules}
\end{figure}

We compare the different ex-p-merging functions introduced in the previous sections, with the addition of the Bonferroni method. The parameter $k$ for the R\"uger combination rule is set to $K/2$ ({twice the median}). The parameter $\rho$ is set to the values $\rho=\{0.1,0.9\}$, corresponding to weak and strong dependence among p-values. Each simulation is repeated for a total of $B=10,000$ replications, and we report the observed empirical average. In Figure~\ref{fig:exch_rules}, we can notice that the error level is controlled at the nominal level $0.05$ for all the proposed methods. Variations in terms of power are observed depending on whether the parameter $\rho$ is set to 0.9 or 0.1. Specifically, {ex-p-merging functions} that exhibit strong performance in the left plot tend to show reduced effectiveness in the right plot, and the opposite is also true. As expected, the Bonferroni method shows a higher power near independence where also the {ex-}Hommel combination rule defined in \eqref{eq:feh*} seems to perform well. A simulation study comparing the different randomized p-merging functions is reported in Appendix \ref{sec:add_sims}.

In this last part, we aim to explore the scenario where p-values are exchangeable under the null hypothesis but not under the alternative. Indeed, the p-values generated as in \eqref{eq:sim_pvals} are exchangeable under both the null and the alternative hypotheses; however, for the results in Section \ref{sec:general_res}, it is only necessary for the p-values to be exchangeable under $H_0$.
Specifically, if the p-values are not exchangeable under the alternative, they could be arranged in a particular way to obtain a more powerful procedure. In other words, if it is suspected that some p-values are smaller under the alternative hypothesis, they could be placed at the beginning of the vector to enhance the power of the procedure since our rules that are valid under exchangeability process the vector of p-values sequentially. 
Clearly, the procedure of ordering p-values in a particular way must preserve the exchangeabilty under the null hypothesis (i.e., data-dependent ordering is usually not allowed since it violates exchangeability).
In the simulation setting, suppose to have $K$ independent studies, each with observations $X_{ij},\,i=1,\dots,K,\,j=1,\dots,n_i,$ that are iid from a normal distribution with mean $\mu$ and variance $1$. In addition, $X_{0j},j=1,\dots,n_0,$ is assumed to be an additional sample from the same population and common for all studies. We define the quantity $\Bar{X}_i$ as
\[
\Bar{X}_i = \frac{1}{\sqrt{n_i}} \sum_{j=1}^{n_i} X_{ij}, \quad i=0,1,\dots,K,
\]
that is distributed as $\mathcal{N}(\sqrt{n_i}\mu,1)$.
The interest is to test the hypothesis $\mu = 0$ under the alternative $\mu \neq 0$
and the test statistic used is 
\[
T_k := \frac{\Bar{X}_k + \Bar{X}_0}{\sqrt{2}}, \quad k=1,\dots,K,
\]
where it is possible to see that each study use the common sample in the ``same way". Under the null hypothesis, the test statistics have the same marginal distribution and the model $(T_1, \dots, T_K)$ is exchangeable. 
Under the alternative, the mean of the test statistics depends on the sample size, so the model cannot be exchangeable. In particular, studies with a higher number of observations are more powerful. The p-values to test the null hypothesis are given by
$$P_k=2\Phi(-|T_k|), \quad k=1,\dots,K.$$

\begin{figure}
    \centering
    \includegraphics[width=0.9\linewidth]{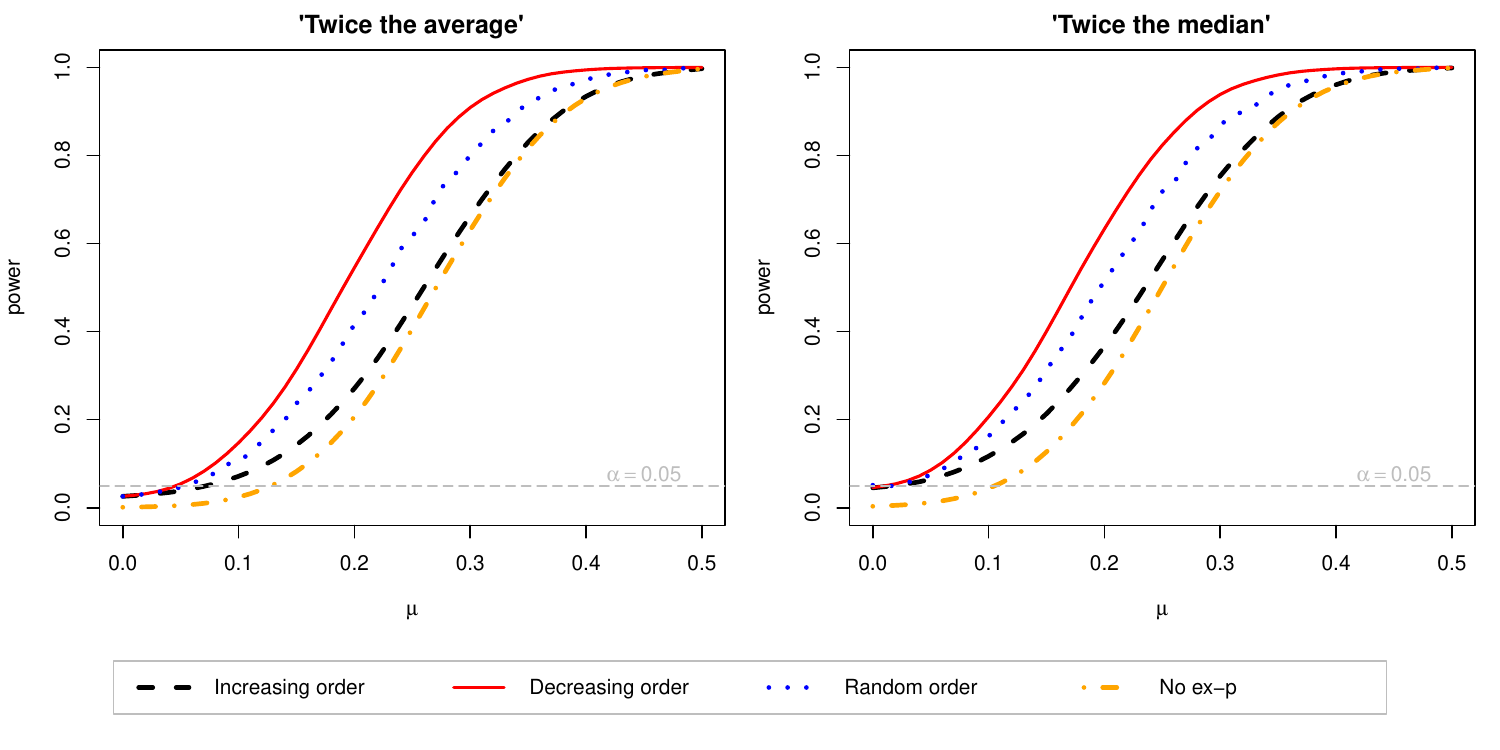}
    \caption{Combination of p-values using different ex-p-merging functions and different ordering based on the sample size. Non ex-p-merging functions valid under arbitrary dependence are added for comparison. The ex-p-merging rules are more powerful if p-values are ordered in decreasing order with respect to the sample size. }
    \label{fig:exch_h0}
\end{figure}

Specifically, in the simulated scenario, $K=10$, $n_i,\,i=1,\dots,K$, take values $10,20,\dots,100$ and $n_0=25$. 
The number of replications is $B=10,000$ and the ex-p-merging functions used are $F_{\mathrm{EA}}$ (``twice the mean") and $F_{\mathrm{ER}}$ with $k=K/2$ (``twice the median"). We let the parameter $\mu$ varies in the interval $[0, 0.5]$ and three different solutions are compared: (i) the p-values are ordered in increasing order with respect to the sample size, (ii) the p-values are ordered in decreasing order with respect to the sample size, and (iii) the p-values are randomly ordered.
In addition, we compare the ex-p-merging functions with the ``standard" rules valid under arbitrary dependence $F_{\mathrm{A}}$ and $F_{\mathrm{R}}$.
The results are reported in Figure \ref{fig:exch_h0}, where we see that when the p-values are ordered in decreasing order with respect to the number of observations, the combined tests are more powerful. This is because the power of individual p-values increases with the sample size. 
Overall, the proposed ex-p-merging are more powerful than the rules valid under arbitrary dependence.

Results under other setups are reported in Appendix \ref{sec:add_sims}, with similar qualitative conclusions.

\section{Summary}\label{sec:summary}
In this paper, we derive novel {p-merging functions} for the scenario where the p-values are exchangeable. These new rules are demonstrated to dominate their original counterparts derived under the assumption of arbitrary dependence. Furthermore, we illustrate how a simple randomization trick (introduction of a uniform random variable or uniform permutation) can also be employed in the case of arbitrarily dependent p-values to yield more powerful rules than existing ones. These results are proposed in a fully general form, addressing the relationship between p-merging functions and e-values introduced by their respective calibrators. In particular, once the corresponding e-value is obtained, it becomes feasible to utilize Markov's inequality or its randomized and exchangeable generalizations from~\cite{ramdas2023}. 

As a practical recommendation, we suggest the exchangeable improvement of the Hommel combination if we have no apriori idea how strong the dependence is likely to be, but to use the exchangeable improvement of ``twice the median'' if the p-values are thought to be strongly dependent.

As an open question, it remains unclear whether our methods are further improvable, stemming from the unknown admissibility or tightness of exchangeable Markov's inequality beyond extreme cases.

\subsection*{Acknowledgments}
The authors thank Ilmun Kim for a suggestion on asymptotic p-values.
\subsection*{Funding}
RW is supported by the
Natural Sciences and Engineering Research Council of Canada (RGPIN-2024-03728, CRC-2022-00141). AR acknowledges NSF grant DMS-2310718 and the Sloan fellowship FG-2024-22012.
\subsection*{Data availability}
The code is available in the GitHub repository: \href{https://github.com/matteogaspa/CombiningExchangeablePValues}{matteogaspa/CombiningExchangeablePValues}.

\bibliography{biblio.bib}

\begin{thebibliography}{}

\bibitem[Banerjee et~al., 2019]{banejeree2019}
Banerjee, M., Durot, C., and Sen, B. (2019).
\newblock {Divide and conquer in nonstandard problems and the super-efficiency phenomenon}.
\newblock {\em The Annals of Statistics}, 47(2):720 -- 757.

\bibitem[Basu, 1955]{basu1955}
Basu, D. (1955).
\newblock On statistics independent of a complete sufficient statistic.
\newblock {\em Sankhyā: The Indian Journal of Statistics (1933-1960)}, 15(4):377--380.

\bibitem[Benjamini and Yekutieli, 2001]{benjamiyekuteli2001}
Benjamini, Y. and Yekutieli, D. (2001).
\newblock The control of the false discovery rate in multiple testing under dependency.
\newblock {\em The Annals of Statistics}, 29(4):1165--1188.

\bibitem[Chen et~al., 2023]{chen2023trade}
Chen, Y., Liu, P., Tan, K.~S., and Wang, R. (2023).
\newblock Trade-off between validity and efficiency of merging p-values under arbitrary dependence.
\newblock {\em Statistica Sinica}, 33(2):851--872.

\bibitem[Chen et~al., 2024]{chen2024sub}
Chen, Y., Wang, R., Wang, Y., and Zhu, W. (2024).
\newblock Sub-uniformity of harmonic mean p-values.
\newblock {\em arXiv preprint arXiv:2405.01368}.

\bibitem[Chi et~al., 2024]{chi2022multiple}
Chi, Z., Ramdas, A., and Wang, R. (2024).
\newblock Multiple testing under negative dependence.
\newblock {\em Bernoulli}, forthcoming.

\bibitem[Choi and Kim, 2023]{choi2023}
Choi, W. and Kim, I. (2023).
\newblock Averaging p-values under exchangeability.
\newblock {\em Statistics \& Probability Letters}, 194:109748.

\bibitem[Cox, 1975]{cox1975}
Cox, D.~R. (1975).
\newblock A note on data-splitting for the evaluation of significance levels.
\newblock {\em Biometrika}, 62(2):441--444.

\bibitem[DiCiccio et~al., 2020]{diciccio2020}
DiCiccio, C.~J., DiCiccio, T.~J., and Romano, J.~P. (2020).
\newblock Exact tests via multiple data splitting.
\newblock {\em Statistics \& Probability Letters}, 166:108865.

\bibitem[Efron, 2010]{efron2012large}
Efron, B. (2010).
\newblock {\em Large-scale inference: empirical Bayes methods for estimation, testing, and prediction}.
\newblock Cambridge University Press.

\bibitem[Fisher, 1934]{fisher1928}
Fisher, R.~A. (1934).
\newblock {\em Statistical methods for research workers}.
\newblock Number~5. Oliver and Boyd.

\bibitem[Grünwald et~al., 2024]{grunwald_safe_2019}
Grünwald, P., de~Heide, R., and Koolen, W. (2024).
\newblock Safe testing.
\newblock {\em Journal of the Royal Statistical Society Series B: Statistical Methodology}, 86(5):1091--1128.

\bibitem[Gui et~al., 2023]{gui2023aggregating}
Gui, L., Jiang, Y., and Wang, J. (2023).
\newblock Aggregating dependent signals with heavy-tailed combination tests.
\newblock {\em arXiv preprint arXiv:2310.20460}.

\bibitem[Guo and Shah, 2024]{guo2023}
Guo, F.~R. and Shah, R.~D. (2024).
\newblock {Rank-transformed subsampling: inference for multiple data splitting and exchangeable p-values}.
\newblock {\em Journal of the Royal Statistical Society Series B: Statistical Methodology}.

\bibitem[Hommel, 1983]{hommel1983}
Hommel, G. (1983).
\newblock Tests of the overall hypothesis for arbitrary dependence structures.
\newblock {\em Biometrical Journal}, 25(5):423--430.

\bibitem[Ignatiadis et~al., 2024]{ignatiadis2024}
Ignatiadis, N., Wang, R., and Ramdas, A. (2024).
\newblock E-values as unnormalized weights in multiple testing.
\newblock {\em Biometrika}, 111(2):417--439.

\bibitem[Kim and Ramdas, 2024]{Kim2024}
Kim, I. and Ramdas, A. (2024).
\newblock {Dimension-agnostic inference using cross U-statistics}.
\newblock {\em Bernoulli}, 30(1):683 -- 711.

\bibitem[Lei and Sudijono, 2024]{lei2024inference}
Lei, L. and Sudijono, T. (2024).
\newblock Inference for synthetic controls via refined placebo tests.
\newblock {\em arXiv preprint arXiv:2401.07152}.

\bibitem[Meinshausen et~al., 2009]{meinshausen2009}
Meinshausen, N., Meier, L., and B{\"u}hlmann, P. (2009).
\newblock P-values for high-dimensional regression.
\newblock {\em Journal of the American Statistical Association}, 104(488):1671--1681.

\bibitem[Morgenstern, 1980]{morgenstern1980}
Morgenstern, D. (1980).
\newblock Berechnung des maximalen signifikanzniveaus des testes "{L}ehne ${H}_0$ ab, wenn $k$ unter $n$ gegebenen tests zur ablehnung f{\"u}hren”.
\newblock {\em Metrika}, 27:285--286.

\bibitem[Owen, 2009]{owen2009}
Owen, A.~B. (2009).
\newblock {Karl Pearson’s meta-analysis revisited}.
\newblock {\em The Annals of Statistics}, 37(6B):3867 -- 3892.

\bibitem[Pearson, 1934]{pearson1934}
Pearson, K. (1934).
\newblock On a new method of determining ``goodness of fit".
\newblock {\em Biometrika}, 26(4):425--442.

\bibitem[Ramdas and Manole, 2024]{ramdas2023}
Ramdas, A. and Manole, T. (2024).
\newblock Randomized and {E}xchangeable {I}mprovements of {M}arkov's, {C}hebyshev's and {C}hernoff's inequalities.
\newblock {\em Statistical Science}, forthcoming.

\bibitem[Ritzwoller and Romano, 2023]{ritzwoller2023}
Ritzwoller, D.~M. and Romano, J.~P. (2023).
\newblock Reproducible aggregation of sample-split statistics.
\newblock {\em arXiv preprint arXiv:2311.14204}.

\bibitem[R{\"u}ger, 1978]{ruger1978}
R{\"u}ger, B. (1978).
\newblock Das maximale signifikanzniveau des tests: “{L}ehne ${H}_0$ ab, wenn $k$ unter $n$ gegebenen tests zur ablehnung f{\"u}hren”.
\newblock {\em Metrika}, 25:171--178.

\bibitem[Rüschendorf, 1982]{ruschendorf1982}
Rüschendorf, L. (1982).
\newblock Random variables with maximum sums.
\newblock {\em Advances in Applied Probability}, 14(3):623--632.

\bibitem[Sarkar, 1998]{sarkar1998}
Sarkar, S.~K. (1998).
\newblock Some probability inequalities for ordered {MTP2} random variables: A proof of the simes conjecture.
\newblock {\em The Annals of Statistics}, 26(2):494--504.

\bibitem[Sellke et~al., 2001]{sellke2001calibration}
Sellke, T., Bayarri, M.~J., and Berger, J.~O. (2001).
\newblock Calibration of $\rho$ values for testing precise null hypotheses.
\newblock {\em The American Statistician}, 55(1):62--71.

\bibitem[Severini, 2000]{severini2000likelihood}
Severini, T.~A. (2000).
\newblock {\em Likelihood methods in statistics}.
\newblock Oxford University Press.

\bibitem[Shafer and Vovk, 2008]{shafer2008}
Shafer, G. and Vovk, V. (2008).
\newblock A {T}utorial on {C}onformal {P}rediction.
\newblock {\em Journal of Machine Learning Research}, 9(3).

\bibitem[Shekhar et~al., 2022]{shekhar2022permutation}
Shekhar, S., Kim, I., and Ramdas, A. (2022).
\newblock A permutation-free kernel two-sample test.
\newblock {\em Advances in Neural Information Processing Systems}, 35:18168--18180.

\bibitem[Shekhar et~al., 2023]{shekhar2023permutation}
Shekhar, S., Kim, I., and Ramdas, A. (2023).
\newblock A permutation-free kernel independence test.
\newblock {\em Journal of Machine Learning Research}, 24(369):1--68.

\bibitem[Simes, 1986]{simes1986}
Simes, R.~J. (1986).
\newblock An improved {B}onferroni procedure for multiple tests of significance.
\newblock {\em Biometrika}, 73(3):751--754.

\bibitem[Stevens, 1950]{stevens1950}
Stevens, W.~L. (1950).
\newblock Fiducial limits of the parameter of a discontinuous distribution.
\newblock {\em Biometrika}, 37(1/2):117--129.

\bibitem[Vovk et~al., 2022]{vovk2022}
Vovk, V., Wang, B., and Wang, R. (2022).
\newblock {Admissible ways of merging p-values under arbitrary dependence}.
\newblock {\em The Annals of Statistics}, 50(1):351--375.

\bibitem[Vovk and Wang, 2020]{vovk2020}
Vovk, V. and Wang, R. (2020).
\newblock Combining p-values via averaging.
\newblock {\em Biometrika}, 107(4):791--808.

\bibitem[Vovk and Wang, 2021]{vovk2021evals}
Vovk, V. and Wang, R. (2021).
\newblock {E-values: Calibration, combination and applications}.
\newblock {\em The Annals of Statistics}, 49(3):1736--1754.

\bibitem[Wang and Wang, 2016]{wang2016joint}
Wang, B. and Wang, R. (2016).
\newblock Joint mixability.
\newblock {\em Mathematics of Operations Research}, 41(3):808--826.

\bibitem[Wang, 2024]{wang2024testing}
Wang, R. (2024).
\newblock Testing with p*-values: Between p-values, mid p-values, and e-values.
\newblock {\em Bernoulli}, 30(2):1313--1346.

\bibitem[Wasserman et~al., 2020]{wasserman2020pnas}
Wasserman, L., Ramdas, A., and Balakrishnan, S. (2020).
\newblock Universal inference.
\newblock {\em Proceedings of the National Academy of Sciences}, 117(29):16880--16890.

\bibitem[Wasserman and Roeder, 2009]{wasserman2009}
Wasserman, L. and Roeder, K. (2009).
\newblock High-dimensional variable selection.
\newblock {\em The Annals of Statistics}, 37(5A):2178--2201.

\bibitem[Westfall et~al., 2004]{westfall2004}
Westfall, P.~H., Kropf, S., and Finos, L. (2004).
\newblock Weighted {FWE}-controlling methods in high-dimensional situations.
\newblock {\em Lecture Notes-Monograph Series}, pages 143--154.

\bibitem[Wilson, 2019]{wilson2019}
Wilson, D.~J. (2019).
\newblock The harmonic mean p-value for combining dependent tests.
\newblock {\em Proceedings of the National Academy of Sciences}, 116(4):1195--1200.

\bibitem[Xu and Ramdas, 2023]{xu2023}
Xu, Z. and Ramdas, A. (2023).
\newblock More powerful multiple testing under dependence via randomization.
\newblock {\em arXiv preprint arXiv:2305.11126}.

\bibitem[Xu et~al., 2024]{xu2022post}
Xu, Z., Wang, R., and Ramdas, A. (2024).
\newblock Post-selection inference for e-value based confidence intervals.
\newblock {\em Electronic Journal of Statistics}, 18(1):2292--2338.

\end{thebibliography}

\newpage
\appendix
\section{Proofs of the results}\label{sec:ap_proof}
\subsection{Proof of Section \ref{sec:notation}}
\begin{proof}[Proof of Lemma~\ref{lemma:eval}]
    By definition, the quantity in~\eqref{eq:sum_of_cals} is non-negative. In addition, for any $\alpha \in (0,1]$,
    \[
    \begin{split}
        \Ev\left[ \frac{1}{\alpha} \sum_{k=1}^K \lambda_k f_k\left(\frac{P_k}{\alpha}\right) \right] &=  \frac{1}{\alpha} \sum_{k=1}^K \lambda_k \Ev\left[ f_k\left(\frac{P_k}{\alpha}\right) \right] = \frac{1}{\alpha} \sum_{k=1}^K \lambda_k \int_0^\alpha f_k\left(\frac{p}{\alpha}\right) \d p\\
        &= \sum_{k=1}^K \lambda_k \int_0^1 f_k\left(p\right) \d p \leq 1.
    \end{split}
    \]
    If the calibrators are admissible, one can  see that the equality holds since $\int_0^1 f_k\left(p\right) \d p=1$ for each $k$.
\end{proof}

\subsection{Proofs of Section \ref{sec:general_res}}
\begin{proof}[Proof of Theorem \ref{th:exch}]
    The proof involves the use of the exchangeable Markov inequality (EMI) recalled in Theorem~\ref{th:EMI} for finite sequences:
    \[
    \pr\left(\exists k \leq K: \frac{1}{k} \sum_{i=1}^k f\left(\frac{P_i}{\alpha}\right) \geq 1 \right) \stackrel{(i)}{\leq} \Ev\left[f\left(\frac{P_1}{\alpha}\right)\right] = \alpha\Ev\left[\frac{1}{\alpha}f\left(\frac{P_1}{\alpha}\right)\right] \stackrel{(ii)}{\leq} \alpha,
    \]
    where $(i)$ is due to EMI while $(ii)$ holds due to Lemma~\ref{lemma:eval}.
\end{proof}

\begin{proof}[Proof of Theorem \ref{th:calibr_to_F_exc}]
    It is clear that $F$ is increasing and Borel since $R_\alpha$ is a lower set. For an exchangeable sequence $\mathbf{P} \in \mathcal{U}^K$ and $\alpha \in (0,1)$, using Theorem~\ref{th:exch} and the fact that $(R_\beta)_{\beta\in(0,1)}$ is nested, we have 
    \begin{equation*}
    \begin{split}
        \pr\left(F(\mathbf{P}) \leq \alpha \right) &= \pr\big( \inf\{\beta \in (0,1): \mathbf{P} \in R_\beta \} \leq \alpha \big)\\
        &= \pr\left(\inf \left\{ \beta \in (0,1): \exists k \leq K \mathrm{~ such~that~} \frac{1}{k} \sum_{i=1}^k f\left(\frac{P_i}{\beta}\right) \geq 1 \right\} \leq \alpha \right)\\
        &= \pr\left( \bigcap_{\beta>\alpha} \left\{ \exists k \leq K: \left(\frac{1}{k} \sum_{i=1}^k f\left( \frac{P_i}{\beta} \right) \right) \geq 1 \right\}  \right) \\
        &= \inf_{\beta > \alpha} \pr \left(\exists k \leq K: \left( \frac{1}{k} \sum_{i=1}^k f\left(\frac{P_i}{\beta} \right) \right) \geq 1 \right) \leq \inf_{\beta > \alpha} \beta = \alpha.
    \end{split}
    \end{equation*}
    Therefore $F$ is a valid p-merging function. Homogeneity comes directly from the definition of~\eqref{eq:f_exch}.
\end{proof}

\begin{proof}[Proof of Proposition \ref{prop:ex-p-m}]
Let $\mathbf P\in \mathcal U^K$,
and let $\sigma$ be a random permutation of $\{1,\dots,K\}$, uniformly drawn from all permutations of $\{1,\dots,K\}$ and independent of $\mathbf P$.
Let $\mathbf P^\sigma= (P_{\sigma(1)}, \dots, P_{\sigma (K)})$. 
Note that $\mathbf P^\sigma$  is exchangeable by construction.
If $F$ is a symmetric ex-p-merging function,
it must satisfy $F(\mathbf P^{\sigma})=F(\mathbf P)$.
Because $F(\mathbf P^{\sigma})$ is a p-variable, 
so is $ F(\mathbf P)$, showing that $F$ is a p-merging function. 
\end{proof}

\begin{proof}[Proof of Proposition \ref{prop:sharp}]
Take  $U\in\mathcal U$ and an event $A$ with $\pr(A)=\alpha$ independent of $U$.
Let $b=f(0+)\le K$. 
The condition on $f$ guarantees that $f(U)$ 
is a random variable  with support $[0,b]$,  mean $1$  and a decreasing density.
The above conditions, using Theorem 3.2 of \cite{wang2016joint}, 
  guarantee that there exists $\mathbf U=(U_1,\dots,U_K)\in \mathcal U^K$ such that 
    $\pr(\sum_{i=1}^K f(U_i) =K)=1$. 
  We   assume that $U,A,\mathbf U$ are mutually independent; this is possible as we are only concerned with distributions.
    Taking a uniformly drawn random permutation further allows us to assume that $\mathbf U$   is exchangeable. 
    Let $P_i=\alpha U_i  \ind_A + (\alpha+(1-\alpha) U)\ind_{A^c} $ for $i=1,\dots,K$. It is clear that each $P_i\in \mathcal U$ and $\mathbf P=(P_1,\dots,P_K)$ is exchangeable. Moreover, by the definition of $F$, $$\pr(F(\mathbf P) \le \alpha) \ge  \pr 
    \left (\frac 1 K\sum_{i=1}^K f(P_i/\alpha  ) \ge 1 \right) = \pr(A)\pr\left (\sum_{i=1}^K f(U_i) =K \right )=\alpha.$$
    The other inequality $\pr(F(\mathbf P) \le \alpha) \le\alpha$ follows from Theorem~\ref{th:exch}.
\end{proof}

\begin{proof}[Proof of Theorem \ref{th:r1-1}]
{We say that a set $R\subseteq [0,\infty)^K$ is a decreasing set if
$\mathbf x \in R $ implies $\mathbf y\in R $ for all
$ \mathbf y\in [0,\infty)^K$ with $\mathbf y
\le  \mathbf x$ (componentwise).

Fix $\alpha\in (0,1) $,
and note that
$R_\alpha(F)= \{\mathbf p \in [0,\infty)^K:  F(\mathbf p)\le \alpha\}$
is a decreasing set. 
Define 
$$G (\mathbf p) = \inf\left\{ \beta \in (0,1) : \mathbf p\in \frac{\beta}{\alpha}  R_\alpha (F)\right \}, ~~~\mathbf p \in [0,\infty)^K. $$
First, $G$ is homogeneous, which follows from its definition. 
Second,   $\mathbf p\in   R_\alpha (F)$ implies $G(\mathbf p)\le \alpha$, and hence $ R_\alpha(F) \subseteq R_\alpha(G)$. 
Third, $G$ is increasing because   $R_\alpha(F)$ is a decreasing set. 

It remains to show that $G$ is a p-merging function. 
The definition of $G$ gives $$G(\mathbf p)\le \beta \iff \mathbf  p \in  \bigcap_{\gamma >\beta} \frac{\gamma }{\alpha} R_\alpha(F). $$
If   we can show
\begin{align}
\label{eq:r1-homo}
\p \left( \mathbf P  \in \frac{\beta}{\alpha}  R_\alpha (F)\right) \le \beta  \mbox{~~~for all  $\mathbf P\in \mathcal U^K$ and $\beta \in (0,1)$},
\end{align}
then
$$
\p ( G(\mathbf P)\le \beta ) \le \inf_{\gamma >\beta} \p \left( \mathbf P  \in  \frac{\gamma }{\alpha}  R_\alpha (F)\right) \le \inf_{\gamma >\beta} \gamma  =\beta,
$$

\begin{lemma}
\label{lem:r1-1}
Let $R\subseteq [0,\infty)^K$ be a decreasing Borel set. 
For any $\beta\in (0,1)$, we have 
$$
\sup_{\mathbf P\in \mathcal U^K} \p(\mathbf P\in \beta R) \ge \beta \iff
\sup_{\mathbf P\in \mathcal U^K} \p(\mathbf P\in   R)  = 1. 
$$
\end{lemma}
\begin{proof}[Proof of the lemma]  
Let $\mathcal V$ be the set of p-variables for $\p$.
For any decreasing set $L$, we have 
\begin{equation}\label{eq:r1-larger}
\sup_{\mathbf P\in \mathcal U^K} \p(\mathbf P\in   L) =\sup_{\mathbf P\in \mathcal V^K} \p(\mathbf P\in   L).
\end{equation}
 This fact will be repeatedly used in the proof below.  

We first prove the $\Leftarrow$ direction  
  by   contraposition. 
Suppose $$\gamma: =\sup_{\mathbf P\in \mathcal U^K}  \p(\mathbf P\in \beta R)< \beta.$$
Take an  event $A$ with probability $\beta$ and any   $\mathbf P \in \mathcal U^K$  independent of $A$. Define $\mathbf P^* $ by $$\mathbf P^*= 
\beta  \mathbf P \times  \ind_A+ \mathbf 1 \times \ind_{A^c}.$$  
It is straightforward to check $\mathbf P^*\in \mathcal V^K$. Hence,  by \eqref{eq:r1-larger},
$$\beta 
\p\left (    \mathbf P\in     R \right)
= \p(A) \p\left (    \mathbf P\in    R  \right)
\le  \p( \mathbf P^*    \in  \beta  R  )  \le \gamma ,$$
and  thus $\p(\mathbf P\in   R)  \le \gamma/\beta$. 
Since  $\gamma/\beta<1$, this yields  $\sup_{\mathbf P\in \mathcal U^K}  \p(\mathbf P\in   R) < 1$  and completes the $\Leftarrow$ direction. 

Next  we show the $\Rightarrow$ direction. 
Suppose   $\sup_{\mathbf P\in \mathcal U^K} \p(\mathbf P\in \beta R) \ge \beta$.
For any $\epsilon\in (0,\beta)$, there exists $\mathbf P=(P_1,\dots,P_K)\in \mathcal U^K$ such that $  \p(\mathbf P\in \beta R) > \beta-\epsilon.$
Let $A=\{\mathbf P\in \beta R\}$, $\gamma=\p(A)$,
and  $B$ be an event containing $A$ with $\p(B)=\beta\vee \gamma$. Let $\mathbf P^*=(P_1^*,\dots,P_K^*)$ follow the conditional distribution of $\mathbf P/\beta$ given $B$. We have
$$
\p( \mathbf P^* \in R) = \p(\mathbf P \in \beta R\mid B) =\p(A \mid B) = \frac{\gamma}{\beta\vee \gamma}.
$$ 
Note that  for $k\in \{1,\dots,K\}$, 
$$
\p(P_k^* \le  p ) = \p(P_k /\beta \le p \mid B) \le \frac{\p(P_k\le \beta p )}{\p(B)} 
=  \frac{\beta p }{\beta\vee \gamma}  \le p,
$$
and hence   $\mathbf P^*\in \mathcal V^K$.
Since $\gamma>\beta-\epsilon $ and $\epsilon\in (0,\beta)$ is arbitrary, we can conclude      
$ 
\sup_{\mathbf P\in \mathcal V^K} \p(\mathbf P\in   R) = 1$, yielding  $\sup_{\mathbf P\in \mathcal U^K} \p(\mathbf P\in   R)  = 1
$   via \eqref{eq:r1-larger}.
\end{proof}
 Now we resume the proof of Theorem \ref{th:r1-1}.  
Fix   $\beta\in (0,1)$.
Take any $\lambda>1$ such that $\lambda(\beta\vee \alpha) <1$. 
 Lemma \ref{lem:r1-1} yields,  for any decreasing set $R$, 
 \begin{align} \label{eq:r1-2}
\sup_{\mathbf P\in \mathcal U^K} \p(\mathbf P\in\lambda \beta R)  < \lambda\beta \iff
\sup_{\mathbf P\in \mathcal U^K} \p(\mathbf P\in \lambda \alpha R) < \lambda \alpha. 
\end{align}
Let $R=R_\alpha(F)/(\alpha \lambda)$, which is a decreasing set. 
Note that 
$$ \sup_{\mathbf P\in \mathcal U^K} \p(\mathbf P\in \lambda \alpha R) = 
 \sup_{\mathbf P\in \mathcal U^K} \p(\mathbf P\in  R_\alpha(F)) \le \alpha < \lambda \alpha 
 $$
 and this leads to, by using \eqref{eq:r1-2},
 $$  \sup_{\mathbf P\in \mathcal U^K} \p\left(\mathbf P\in \frac{\beta}{\alpha} R_\alpha(F)\right)  =
 \sup_{\mathbf P\in \mathcal U^K} \p(\mathbf P\in\lambda \beta R)  < \lambda\beta. 
 $$
 Since $\lambda>1$ can be arbitrarily close to $1$, we conclude that \eqref{eq:r1-homo} holds true,
 and this completes the proof.}
     \end{proof}

\begin{proof}[Proof of Lemma \ref{lemma:continuous}]
    Fix any $\mathbf p = (\mathbf p_1, \mathbf p_2 ) \in [0,1]^K$ and suppose that $F(\mathbf{p}_1) = \alpha$. By definition, we have
    \[
    \exists k \leq K_1: \frac 1k \sum_{i=1}^k f\left(\frac{p_i}{\alpha} \right) \geq 1 \implies \exists k \leq K_1 + K_2: \frac1k \sum_{i=1}^k f\left(\frac{p_i}{\alpha} \right) \geq 1.
    \]
    This implies 
    \[
    F(\mathbf p) = \inf\left\{\beta \in (0,1): \exists k \leq K \mathrm{~s.t.~}\frac{1}{k} \sum_{i=1}^k f\left(\frac{p_i}{\beta} \right) \geq 1 \right\} \leq \alpha,
    \]
    due to the fact that $f$ is increasing in $\beta$.
\end{proof}

\begin{proof}[Proof of Theorem \ref{th:sequential}]
    From Lemma~\ref{lemma:continuous}, we have that 
    \[
    \exists k \leq K: F(\mathbf{p}_k) \leq \alpha \iff F(\mathbf{p}_K) \leq \alpha,
    \]
    where $\mathbf{p}_k = (p_1,\dots,p_k)$ is the sequence containing the first $k$ values of $\mathbf{p}$. Then we can write,
    \[
    \begin{split}
    \pr\left(\exists k \geq 1: F(\mathbf{P}_k) \leq \alpha\right) &= \pr\left(\bigcup_{K \in \mathbb{N}} \{ \exists k \leq K: F(\mathbf{P}_k) \leq \alpha\}\right) \\
    &= \lim_{K \to \infty} \pr\left(\exists k \leq K: F(\mathbf{P}_k) \leq \alpha\right)\\
    &= \lim_{K \to \infty} \pr\left(F(\mathbf{P}_K) \leq \alpha\right) \leq \alpha,
    \end{split}
    \]
    where the last inequality is due to Theorem~\ref{th:calibr_to_F_exc}.
\end{proof}

\begin{proof}[Proof of Theorem \ref{th:1}]
From direct calculation and using Lemma~\ref{lemma:eval},
\begin{align*}
\pr   
    \left(
      \sum_{k=1}^K \lambda_k f_k\left(\frac{P_k}{\alpha}\right) \ge U
    \right)
    &=\Ev\left[ \pr   
    \left(
      \sum_{k=1}^K \lambda_k f_k\left(\frac{P_k}{\alpha}\right) \ge U
    \right)\mid \mathbf P\right]
    \\&= \Ev \left[  \left(\sum_{k=1}^K \lambda_k f_k\left(\frac{P_k}{\alpha}\right)\right)\wedge 1\right ]
          \\& \le    \Ev \left[  \sum_{k=1}^K \lambda_k f_k\left(\frac{P_k}{\alpha}\right)\right]  = \alpha \Ev\left[\frac{1}{\alpha}\sum_{k=1}^K \lambda_k f_k\left(\frac{P_k}{\alpha} \right)\right] \le \alpha,
\end{align*}
The equality for $\beta =\alpha$ follows because $\sum_{k=1}^K \lambda_k f_k(P_k/\alpha) \le 1$ and $\int_0^1 f_k(p)\d p=1$ for each $k$ guarantee the inequalities in the above set of equations are equalities. For $\beta <\alpha$, it suffices to notice that 
$\sum_{k=1}^K \lambda_k f_k(P_k/\alpha)$ is increasing in $\alpha$.
\end{proof}

\begin{proof}[Proof of Theorem \ref{th:2}]
It is clear that $F$ is increasing and Borel since $R_\alpha$ is a lower set.
For $ (\mathbf P,U)=(P_1,\dots,P_K,U) \in \mathcal U^K\otimes \mathcal U$ and $\alpha \in (0,1)$, using Theorem~\ref{th:1} and the fact that $(R_\beta)_{\beta \in (0,1)}$ is nested,
we have 
\begin{align}
\pr(F(\mathbf P, U) \le\alpha ) &=
\pr\left(\inf\left\{\beta \in (0,1):   \sum_{k=1}^K \lambda_k f_k\left(\frac{P_k}{\beta}\right) \ge U \right\}\le \alpha\right)  \notag
\\&=
\pr\left( \bigcap_{\beta>\alpha} \left\{  \sum_{k=1}^K \lambda_k f_k\left(\frac{P_k}{\beta }\right) \ge U  \right\}\right) \notag
\\&=\inf_{\beta >\alpha} \pr\left(    \sum_{k=1}^K \lambda_k f_k\left(\frac{P_k}{\beta }\right) \ge U   \right)
\le
\inf_{\beta >\alpha} \beta =\alpha. 
\notag 
\end{align}
Therefore, $F$ is a randomized p-merging function.
Homogeneity of $F$ follows  from \eqref{eq:def-f}. 

Since $F$ is homogeneous and increasing, it is continuous in $\mathbf p$. 
Moreover, for fixed $\mathbf p\in [0,1]^K$, since
$$
\bigcap_{v<u} \left\{\alpha \in (0,1):   \sum_{k=1}^K \lambda_k f_k\left(\frac{p_k}{\alpha}\right) \ge v \right\} 
= \left\{\alpha \in (0,1):   \sum_{k=1}^K \lambda_k f_k\left(\frac{p_k}{\alpha}\right) \ge u \right\},
$$
we have 
\begin{align*}
\lim_{v\uparrow u} F(\mathbf p,v) &= \lim_{v\uparrow u} \inf \left\{\alpha \in (0,1):   \sum_{k=1}^K \lambda_k f_k\left(\frac{p_k}{\alpha}\right) \ge v \right\} 
\\&= \inf \bigcap_{v<u}\left\{\alpha \in (0,1):   \sum_{k=1}^K \lambda_k f_k\left(\frac{p_k}{\alpha}\right) \ge v \right\}  
 = F( \mathbf p, u).
\end{align*}
Therefore, $u\mapsto F(\mathbf p,u)$ is lower semi-continuous.   
\end{proof}

\subsection{Proofs of Section \ref{sec:ruger}}

\begin{proof}[Proof of Theorem \ref{thm:ruger_exch}]
According to Theorem~\ref{th:calibr_to_F_exc}, it follows that the function $F_{\mathrm{ER}}$ is a valid ex-p-merging function. 
Fix any $\alpha \in (0,1)$ and $\mathbf{p}\in (0,1]^{K}$. Note that $F_{\mathrm{ER}}(\mathbf{p}) \leq \alpha$ if and only if 
    \[
    \exists \ell \leq K: \frac{1}{\ell} \sum_{i=1}^\ell \frac{K}{k} \ind\left\{\frac{p_i}{\alpha} \leq \frac{k}{K} \right\} \geq 1 \implies \exists \ell \leq K: \sum_{i=1}^\ell \ind\left\{p_i \leq \alpha\frac{k}{K} \right\} \geq \left\lceil \ell \frac{k}{K} \right\rceil.
    \]
    Where rounding up is due to the fact that the summation takes values in positive integers. This holds true if and only if  
    \[
    \exists \ell \leq K: p^\ell_{(\lambda_\ell)} \leq \alpha \frac{k}{K},
    \]
    where $p^\ell_{(\lambda_\ell)}$ is the $\lceil \ell \frac{k}{K} \rceil$-th ordered value of $(p_1, \dots, p_\ell)$. Rearranging the terms, we obtain that it is verified when
    \[
    \frac{K}{k}\bigwedge_{\ell=1}^K p_{(\lambda_\ell)}^\ell \leq \alpha,
    \]
    which complete the first part of the proof. For the second statement, it is possible to note that the element $p_{(\lambda_\ell)}^\ell$ in the sequence coincides with $p_{(k)}$ when $\ell=K$. 
\end{proof}

\begin{proof}[Proof of Theorem \ref{th:ruger_rand}]
According to Corollary~\ref{co:1}, it follows that the function $F_{\mathrm{UR}}$ is a valid randomized p-merging function. 
Fix any $\alpha \in (0,1)$ and $(\mathbf{p}, u) \in (0,1]^{K+1}$, then it is possible to note that $F(\mathbf{p}, u) \leq \alpha$ if and only if
    \[
    \frac{1}{K} \sum_{i=1}^K \frac{K}{k} \ind\left\{\frac{p_i}{\alpha} \leq \frac{k}{K} \right\} \geq u \implies  \sum_{i=1}^K \ind\left\{p_i \leq \alpha \frac{k}{K} \right\} \geq \lceil u k \rceil.
    \]
    Rearranging the terms this holds true only if 
    \[
    p_{(\lceil uk \rceil)} \leq \alpha \frac{k}{K} \implies \frac{K}{k} p_{(\lceil uk \rceil)} \leq \alpha,
    \]
    which concludes the claim. Since $u\leq 1$ almost surely, we have $ p_{(\lceil uk \rceil)} 
 \leq p_{(k)}$. 
\end{proof}

\subsection{Proofs of Section \ref{sec:hommel}}
\begin{proof}[Proof of Lemma~\ref{lemma:calib_hommel}]
    It is simple to see that $f$ is decreasing and it is upper semicontinuous (the function $f$ has discontinuity points in $i/(Kh_K),~i=1,\dots,K,$ but it is simple to prove that $\lim_{x\to x_0} f(x) \leq f(x_0)$). In addition,
    \[
    \int_0^1 f(p) \d p = \int_0^1 \frac{K \ind\{h_K p \leq 1\}}{\lceil K h_K p\rceil} \d p = \sum_{j=1}^K \frac{K}{j} \frac{1}{Kh_K} = \frac{1}{h_K} \sum_{j=1}^K\frac{1}{j} = 1.
    \]
    Due to Theorem~\ref{th:vovkwang} we have that $F_{\mathrm{Hom}}$ is admissible.

    To prove the last part, we see that for any $\mathbf{p} \in (0,1]^K$ and $\alpha \in (0, 1)$ we have that $F_\mathrm{Hom}(\mathbf{p}) \leq \alpha$, if and only if
    \[
    \bigwedge_{k=1}^K \frac{K}{k}p_{(k)} \leq \frac{\alpha}{h_K}.
    \]
    This implies that, for some $m \in \{1, \dots, K\}$, we have that 
    \[
    \sum_{j=1}^K \ind\left\{\frac{K}{m} h_K p_j \leq \alpha \right\} \geq m.
    \]
    We now define this chain of inequalities,
    \[
    1 \leq \sum_{j=1}^K \frac{1}{m} \ind\left\{\frac{K}{m} h_K p_j \leq \alpha \right\} \stackrel{(i)}{\leq} \sum_{j=1}^K \frac{1}{\lceil K h_K p_j / \alpha \rceil} \ind\left\{ \frac{K}{m} h_K p_j \leq \alpha \right\} \stackrel{(ii)}{\leq} \sum_{j=1}^K \frac{1}{\lceil K h_K p_j / \alpha \rceil} \ind\left\{h_K p_j \leq \alpha \right\},
    \]
    where $(i)$ is a consequence of $(1/m)\ind\{K p_j h_K \leq \alpha m\} \leq (1/\lceil K h_K p_j / \alpha \rceil)\ind\{K p_j h_K \leq \alpha m\}$, for all $j=1,\dots,K$, while $(ii)$ to the fact that $K/m \geq 1$. This concludes the proof.
\end{proof}    

\begin{proof}[Proof of Theorem \ref{thm:hom_exch}]
    According to Theorem~\ref{th:calibr_to_F_exc} and Lemma~\ref{lemma:calib_hommel}, it follows that $F_{\mathrm{EHom}}$ is a ex-p-merging function. It is simple to see that $F_{\mathrm{EHom}} \leq F_{\mathrm{Hom}}$.
\end{proof}

\begin{proof}[Proof of Theorem \ref{thm:hom_rand}]
    According to Corollary~\ref{co:1} and Lemma~\ref{lemma:calib_hommel}, it follows that $F_{\mathrm{UHom}}$ is a randomized p-merging function. It is simple to see that $F_{\mathrm{UHom}} \leq F_{\mathrm{Hom}}$.
\end{proof}

\subsection{Proofs of Section \ref{sec:avg}}
\begin{proof}[Proof of Lemma \ref{lemma:calibrator}]
It is easy to see that $f$ is decreasing and
$
\int_0^1 (2-2p) \d p = 1.
$
In addition,  it is upper semi-continuous and $f(0)=\infty$.
\end{proof}

The statements on domination in Theorems \ref{thm:avg_exch}, \ref{theorem:FUA}, and other similar results are straightforward from definitions and we omit the proof.   
\begin{proof}[Proof of Theorem \ref{thm:avg_exch}]
According to Theorem~\ref{th:calibr_to_F_exc} and Lemma~\ref{lemma:calibrator}, it follows that the function $F_{\mathrm{EA}}$ is an ex-p-merging function.
    Fix any $\alpha \in (0,1)$ and $\mathbf{p} \in (0,1]^K$. We note that $F_{\mathrm{EA}}(\mathbf{p}) \leq \alpha$ if and only if, for some $\ell \in \{1, \dots, K\}$, we have
    \begin{equation}\label{eq:th_e_a}
        \frac{1}{\ell} \sum_{i=1}^\ell \left(2 - \frac{2p_i}{\alpha} \right)_+ \geq 1.
    \end{equation}
    This implies that there exists $\ell \leq K$ such that
    \[
    \frac{1}{\ell} \sum_{j=1}^m \left(2 - \frac{2p^\ell_{(j)}}{\alpha} \right) \geq 1 \mathrm{~~~for~some~}m \in \{1, \dots, \ell\},
    \]
    where we recall that $p^\ell_{(j)}$ is the $j$-th ordered value of the vector $(p_1, \dots,p_\ell)$. This is due to the fact that the contribution of $p_i$ in the left-hand side of~\eqref{eq:th_e_a} vanishes for large values of $p_i$. Rearranging the terms, it is possible to obtain that exists $\ell \leq K$ such that
    \[
    \frac{2 A(\mathbf{p}^\ell_{(m)})}{2 - \ell/m}\leq \alpha \mathrm{~~~for~some~}m \in \{1, \dots, \ell\}.
    \]
    Taking an infimum over $m$ yields
    \[
    \left( \bigwedge_{m=1}^\ell \frac{2A(\mathbf{p}^\ell_{(m)})}{(2-\ell/m)_+} \right)\leq \alpha \mathrm{~~~for~some}~\ell\in\{1, \dots, K\}.
    \]
    Actually, the index $m$ can start from $\lceil \ell/2 \rceil$ since the first $\lceil \ell/2 \rceil - 1$ terms in the denominator are smaller than zero. Taking an infimum also over $\ell$ gives the desired result. 
\end{proof}

\begin{proof}[Proof of Theorem \ref{theorem:FUA}]
According to Corollary~\ref{co:1} and Lemma~\ref{lemma:calibrator}, it follows that the function $F_{\mathrm{UA}}$ is a randomized p-merging function.
Fix any $\alpha \in (0,1)$ and $(\mathbf p,u)\in (0,1]^{K+1}$.
Note that $F_{\mathrm{UA}}(\mathbf p,u)\le \alpha$
if and only if $$\frac{1}{K}\sum_{k=1}^m \left(2-\frac{2p_{(k)}}{\alpha}\right) \ge u \mbox{~~~for some $m\in\{1,\dots,K\}$}.$$
Rearranging terms, it is 
$$
\frac{\sum_{k=1}^m {2p_{(k)}}}{2m-Ku}  \le \alpha \mbox{~~~for some $m\in\{1,\dots,K\}$}.
$$
Taking an infimum over $m$ yields the desired formula.
\end{proof}

\subsection{Proofs of Section \ref{sec:harm_avg}}
\begin{proof}[Proof of Lemma~\ref{lemma:calib_harmonic}]
    Let $K \geq 2$ and $T_K \geq 1$. Define the function $f:[0,\infty) \to [0, \infty]$ as
    \[
    p \mapsto \min\left\{\frac{1}{T_Kp} - \frac{1}{T_K}, K\right\} \ind\{p \in [0,1]\},
    \]
that is decreasing in $p$. Then,
\[
\begin{split}
    \int_0^1 f(p) \d p &= \int_0^1 \min\left\{\frac{1}{T_Kp} - \frac{1}{T_K}, K\right\} \ind\{p \in [0,1]\} \d p\\
    &= \int_0^{(T_KK+1)^{-1}} K \d p + \int_{(T_KK+1)^{-1}}^1 \left(\frac{1}{T_Kp} - \frac{1}{T_K} \right) \d p \\
    &= \frac{K}{T_KK+1} - \frac{K}{T_KK+1} + \frac{\log(T_KK+1)}{T_K} = \frac{\log(KT_K+1)}{T_K}.
\end{split}
\]
This implies that $\int_{0}^1 f(p) \d p \leq 1$ if and only if
\begin{equation}\label{eq:proof_c}
     KT_K+1-e^T_K \leq 0.
\end{equation}
We would like to choose $T_K$ as small as possible. One possible candidate is $T_K= \log K + \log \log K +1$. Indeed, plugging $T_K=\log K + \log \log K +1$ into the left-hand side of~\eqref{eq:proof_c} we find that it is verified when
\[
\log \log K + 1 + \frac{1}{K} \leq (e-1) \log K,
\]
and this holds if $K \geq 2$ by checking $K = 2, 3$ and using the derivative of both
sides for $K\geq 4$.
\end{proof}

\begin{proof}[Proof of Proposition~\ref{prop:harm_avg}]
It is straightforward to see that $F$ is an increasing function. By direct calculation,
\begin{equation*}
\begin{split}
    \pr(F(\mathbf{P}) \leq \alpha) &= \pr((T_K+1)H(\mathbf{P}) \leq \alpha)\\
    &= \pr\left( (T_K+1) K \left(\sum_{k=1}^K \frac{1}{P_k} \right)^{-1} \leq \alpha \right)\\
    &= \pr\left(\frac{1}{K} \sum_{k=1}^K \left(\frac{\alpha}{T_K P_k} - \frac{1}{T_K} \right) \geq 1 \right)\\
    &\stackrel{(i)}{\leq} \pr\left(\frac{1}{K} \sum_{k=1}^K \left(\frac{1}{T_K P_k/\alpha} - \frac{1}{T_K} \right)\ind\left\{\frac{P_k}{\alpha} \in [0,1] \right\} \geq 1 \right)\\
    &= \pr\left(\frac{1}{K} \sum_{k=1}^K \min\left\{\left(\frac{1}{T_K P_k/\alpha} - \frac{1}{T_K} \right), K\right\}\ind\left\{\frac{P_k}{\alpha} \in [0,1] \right\} \geq 1 \right)\\
    &\stackrel{(ii)}{\leq} \Ev\left[\frac{1}{K} \sum_{k=1}^K \min\left\{\left(\frac{1}{T_K P_k/\alpha} - \frac{1}{T_K} \right), K\right\}\ind\left\{\frac{P_k}{\alpha} \in [0,1] \right\}\right]\\
    &= \alpha\Ev\left[\frac{1}{\alpha}\frac{1}{K} \sum_{k=1}^K \min\left\{\left(\frac{1}{T_K P_k/\alpha} - \frac{1}{T_K} \right), K\right\}\ind\left\{\frac{P_k}{\alpha} \in [0,1] \right\}\right] \leq \alpha,
\end{split}
\end{equation*}
where $(i)$ is due to the fact that $(1/(T_K x) - 1/T_K)$ is negative for $x>1$, and $(ii)$ holds due to Markov's inequality. The last inequality is a consequence of Lemma~\ref{lemma:eval} and Lemma~\ref{lemma:calib_harmonic}.
\end{proof}

\begin{proof}[Proof of Theorem \ref{thm:harm_exc}]
According to Theorem~\ref{th:calibr_to_F_exc} and Lemma~\ref{lemma:calib_harmonic}, it follows that $F_{\mathrm{EH}}$ is an ex-p-merging function.
Fix any $\alpha \in (0,1)$ and $\mathbf{p} \in (0,1]^K$. We note that $F_{\mathrm{EH}}(\mathbf{p}) \leq \alpha$ if and only if, for some $\ell \in \{1, \dots, K\}$, we have
\[
\frac{1}{\ell} \sum_{i=1}^\ell \left(\frac{\alpha}{T_Kp_i} - \frac{1}{T_K}\right)_+ \geq 1.
\]
This implies that exists $\ell \leq K$ such that
\[
\frac{1}{\ell} \sum_{j=1}^m \left(\frac{\alpha}{T_K p^\ell_{(j)}} - \frac{1}{T_K}\right) \geq 1 \quad \mathrm{~for~some~}m\in\{1, \dots, \ell\},
\]
where we recall that $p^\ell_{(j)}$ is the $j$-th ordered value of the vector $(p_1, \dots, p_\ell)$.
Rearranging the terms it is possible to obtain that exist $\ell \leq K$ such that
\[
\left(\frac{\ell\, T_K}{m} + 1 \right) H(\mathbf{p}^\ell_{(m)}) \leq \alpha \quad \mathrm{~for~some~}m\in\{1, \dots, \ell\},
\]
Taking an infimum over $m$, we get
\[
\bigwedge_{m=1}^\ell \left(\frac{\ell \, T_K}{m} + 1\right) H(\mathbf{p}^\ell_{(m)}) \leq \alpha \quad \mathrm{~for~some~}\ell\in\{1, \dots, K\}.
\]
Taking an infimum over $\ell$ yields the desired result.
\end{proof}

\begin{proof}[Proof of Thereom \ref{thm:harm_rand}]
    According to Corollary~\ref{co:1} and Lemma~\ref{lemma:calib_harmonic}, it follows that $F_{\mathrm{UH}}$ is a randomized p-merging function.
    Fix any $\alpha \in (0,1)$ and $(\mathbf{p}, u) \in (0,1]^{K+1}$. Then $F_{\mathrm{UH}}(\mathbf{p}, u) \leq \alpha$ if and only if
    \[
    \frac{1}{K} \sum_{k=1}^m \left(\frac{\alpha}{T_K p_{(k)}} - \frac{1}{T_K}\right) \geq u \quad \mathrm{~for~some~} m\in\{1, \dots, K\}.
    \]
    Rearranging the terms, it is 
    \[
    \left({uKT_K} + m\right) \left(\sum_{k=1}^m \frac{1}{p_{(k)}}\right)^{-1} \leq \alpha \quad \mathrm{for~some~} m\in\{1, \dots, K\}.
    \]
    Taking a minimum over $m$ yields the desired formula.
\end{proof}

\subsection{Proofs of Section \ref{sec:geo_avg}}
\begin{proof}[Proof of Theorem \ref{thm:geom_exch}]
According to Theorem~\ref{th:calibr_to_F_exc}, it follows that the function $F_{\mathrm{EG}}$ is an ex-p-merging function.
 Fix any $\alpha \in (0,1)$ and $\mathbf{p} \in (0,1]^K$. Then $F_{\mathrm{EG}}(\mathbf{p}) \leq \alpha$, if and only if exists $\ell \in \{1, \dots, K\}$ such that
    \[
    \frac{1}{\ell} \sum_{i=1}^\ell(-\log p_i + \log \alpha)_+ \geq 1.
    \]
    This is verified when exists $\ell \leq K$ such that
    \[
    \frac{1}{\ell} \sum_{j=1}^m(-\log p^\ell_{(j)} + \log \alpha) \geq 1 \quad \mathrm{~for~some~}m \in \{1, \dots, \ell\},
    \]
    where we recall that $p^\ell_{(j)}$ is the $j$-th ordered value of the vector $(p_1, \dots, p_\ell).$ Rearranging the terms it is possible to obtain that exists $\ell\leq K$ such that
    \[
    e^{\ell/m} G(\mathbf{p}^\ell_{(m)}) \leq \alpha \quad \mathrm{~for~some~}m \in \{1, \dots, \ell\}.
    \]
    Taking an infimum over $m$, we get
    \[
    \bigwedge_{m=1}^\ell \left(e^{\ell/m} G(\mathbf{p}^\ell_{(m)}) \right) \leq \alpha \quad \mathrm{~for~some~} \ell \in \{1, \dots, K\}.
    \]
    Taking an infimum over $\ell$ yields the desired result.
\end{proof}

\begin{proof}[Proof of Theorem \ref{thm:geom_rand}]
    According to Corollary~\ref{co:1}, it follows that $F_{\mathrm{UG}}$ is a randomized p-merging function.
    Fix any $\alpha \in (0,1)$ and $(\mathbf{p}, u) \in (0,1]^{K+1}$. We have that $F_{\mathrm{UG}}(\mathbf{p}, u) \leq \alpha$ if and only if
    \[
    \frac{1}{K} \sum_{k=1}^m (-\log p_{(k)} + \log \alpha) \geq u \quad \mathrm{for~some~} m\in\{1, \dots, K\}.
    \]
    Rearranging the terms, it is 
    \[
    e^{u \frac{K}{m}} G(\mathbf{p}_{(m)}) \leq \alpha \quad \mathrm{for~some~} m\in\{1, \dots, K\}.
    \]
    Taking a minimum over $m$ yields the desired formula.
\end{proof}

\section{Exchangeable and randomized p-merging function}\label{sec:exch_rand_pmerg}
In this part, we will integrate the results in Section~\ref{sec:general_res}, using both exchangeability and randomization. In fact, starting from exchangeable p-values, it is possible to prove that if randomization is allowed, then it is possible to improve some of the results obtained in Section \ref{sec:general_res}. We start by defining a ``randomized'' version of Theorem~\ref{th:exch}.
\begin{theorem}\label{th:exch+rand}
    Let $f$ be a calibrator, and $(\mathbf{P}, U) = (P_1, \dots, P_K, U) \in \mathcal{U}^K \otimes \mathcal{U}$ such that $\mathbf{P}$ is exchangeable. For each $\alpha \in (0,1)$, we have
    \[
    \pr\left(f\left(\frac{P_1}{\alpha}\right) \geq U \mathrm{~or~} \exists k \leq K: \frac{1}{k} \sum_{i=1}^k f\left(\frac{P_i}{\alpha}\right) \geq 1  \right) \leq \alpha.
    \]
\end{theorem}
\begin{proof}
    The proof invokes the exchangeable and uniformly-randomized Markov inequality (EUMI) introduced in~\cite{ramdas2023}; see Theorem~\ref{th:EUMI}. In particular,
    \[
    \begin{split}
       \pr\left(f\left(\frac{P_1}{\alpha}\right) \geq U \mathrm{~or~} \exists k \leq K: \frac{1}{k} \sum_{i=1}^k f\left(\frac{P_i}{\alpha}\right) \geq 1 \right) \stackrel{(i)}{\leq} \Ev\left[f\left(\frac{P_1}{\alpha} \right) \right] = \alpha \Ev\left[\frac{1}{\alpha}f\left(\frac{P_1}{\alpha} \right) \right] \stackrel{(ii)}{\leq} \alpha,
    \end{split}
    \]
    where $(i)$ is due to EUMI and $(ii)$ is due to Lemma~\ref{lemma:eval}.
\end{proof}

Similarly to how it was done in the preceding sections, let us now define a randomized ex-p-merging function.

\begin{definition}
    A randomized ex-p-merging function is an increasing Borel function $F:[0,1]^{K+1} \to [0,1]$ such that $\pr(F(\mathbf{P}, U) \leq \alpha) \leq \alpha$ for all $\alpha \in (0,1)$ and $(\mathbf{P}, U) \in \mathcal{U}^K \otimes \mathcal{U}$ with $\mathbf{P}$ exchangeable. It is homogeneous if $F(\gamma \mathbf{p}, u) = \gamma F(\mathbf{p}, u)$ for all $\gamma \in (0,1]$ and $(\mathbf{p}, u) \in [0, 1]^{K+1}$. A randomized ex-p-merging function is admissible if for any randomized ex-p-merging function $G$, $G \leq F$ implies $G=F$.
\end{definition}

 Let $f$ be a calibrator; then for $\alpha \in (0,1)$, we define the exchangeable and randomized rejection region by
 \[
 R_\alpha = \left\{(\mathbf{p},u) \in [0,1]^{K+1}: f\left(\frac{p_1}{\alpha}\right) \geq u \mathrm{~or~} \exists k \leq K: \frac{1}{k} \sum_{i=1}^k f\left(\frac{p_i}{\alpha}\right) \geq 1 \right\},
 \]
 where we set $f(p_i/u)=0$ if $u=0$. Starting from $R_\alpha$, we can define the function $F:[0,1]^{K+1} \to [0,1]$ by
 \begin{equation}\label{eq:exch_rand_f}
 \begin{split}
      F(\mathbf{p}, u) &= \inf\left\{\alpha \in (0,1): (\mathbf{p}, u) \in R_\alpha \right\}\\
      &= \inf\left\{\alpha \in (0,1): f\left(\frac{p_1}{\alpha} \right) \geq u \mathrm{~or~} \exists k \leq K: \frac{1}{k} \sum_{i=1}^k f\left(\frac{p_i}{\alpha} \right) \geq 1 \right\},
 \end{split}
 \end{equation}
with the convention $\inf \varnothing = 1$ and $0 \times \infty = \infty$.
\begin{theorem}
    If $f$ is a calibrator and $(\mathbf{P}, U) \in \mathcal{U}^K \otimes \mathcal{U}$ with $\mathbf{P}$ exchangeable, then $F$ in \eqref{eq:exch_rand_f} is a homogeneous randomized ex-p-merging function.
\end{theorem}
\begin{proof}
    It is clear that $F$ is increasing and Borel since $R_\alpha$ is a lower set. For an exchangeable $\mathbf{P} \in \mathcal{U}^K$ and $\alpha \in (0,1)$, using Theorem~\ref{th:exch+rand} and the fact that $(R_\beta)_\beta \in (0,1)$ is nested, we have
    \[
    \begin{split}
        \pr(F(\mathbf{P}, U) \leq \alpha) &= \pr\left(\inf\left\{\beta \in (0,1): f\left( \frac{P_1}{\beta}\right) \geq U \mathrm{~or~} \exists k \leq K: \frac{1}{k} \sum_{i=1}^k f\left(\frac{P_i}{\beta} \right) \geq 1 \right\} \leq \alpha \right)\\
        &= \pr\left(\bigcap_{\beta > \alpha} \left\{f\left( \frac{P_1}{\beta}\right) \geq U \mathrm{~or~} \exists k \leq K: \frac{1}{k} \sum_{i=1}^k f\left(\frac{P_i}{\beta} \right) \geq 1 \right\} \right)\\
        &= \inf_{\beta > \alpha} \pr\left(f\left( \frac{P_1}{\beta}\right) \geq U \mathrm{~or~} \exists k \leq K: \frac{1}{k} \sum_{i=1}^k f\left(\frac{P_i}{\beta} \right) \geq 1\right)\\
        &\leq \inf_{\beta > \alpha} \beta = \alpha.
    \end{split}
    \]
    therefore $F$ is a valid randomized ex-p-merging function. Homogeneity comes directly from the definition of~\eqref{eq:exch_rand_f}.
\end{proof}

\section{Generalized Hommel combination rule}
\label{sec:app_genHom}
We can generalize the Hommel combination by allowing the selection of certain quantiles from the possible $K$ different quantiles of $\mathbf{p}=(p_1, \dots, p_K)$, for example, one can select the minimum between $K$ times the minimum, 2 times the median and the maximum. This can be obtained by noting that the Hommel combination rule can be rewritten in terms of quantiles. In particular, the following holds:
\[
F'_{\mathrm{Hom}}(\mathbf{p}) = h_K  \bigwedge_{k=1}^K \frac{1}{\lambda_k} p_{(\lceil \lambda_k K\rceil)}, ~~~~\mbox{with }h_K= \sum_{j=1}^K \frac{\lambda_{j} - \lambda_{j-1}}{\lambda_j},
\]
where $(\lambda_0, \lambda_1, \dots, \lambda_K)$ is the vector of quantiles such that $\lambda_j = j/K, \, j=0,1,\dots,K$. This gives the intuition to define a generalization of the aforementioned rule.

Let us define the vector of quantiles $\lambda=(\lambda_0, \lambda_1, \dots, \lambda_M)$, such that $\lambda_0=0,\,\lambda_j\in(0,1],$ 
 if $j=1,\dots,M$ and $\lambda_j < \lambda_{j+1}$. Then, we can define
\begin{equation}\label{eq:ghom}
    F'_{\mathrm{GHom}}(\mathbf{p}) := h_M \bigwedge_{k=1}^M \frac{1}{\lambda_k} p_{(\lceil \lambda_k K \rceil)}, ~~~~\mbox{with }h_M= \sum_{j=1}^M \frac{\lambda_{j} - \lambda_{j-1}}{\lambda_j}.
\end{equation}

\begin{lemma}
    Let $f$ be a function defined by    \begin{equation*}
        f(p) = \sum_{j=1}^M \frac{1}{\lambda_j}\ind\left\{p \in \left(\frac{\lambda_{j-1}}{h_M},  \frac{\lambda_{j}}{h_M}\right] \right\} + \infty \ind\{ p = 0\}.
    \end{equation*}
    Then $f$ is an admissible calibrator. Moreover, the p-merging function induced by $f$ is
    \begin{equation*}
            F_\mathrm{GHom}(\mathbf{p})= \inf\left\{\alpha \in (0,1): \frac{1}{K} \sum_{k=1}^K \sum_{j=1}^M \frac{1}{\lambda_j}\ind\left\{\frac{p_k}{\alpha} \in \left(\frac{\lambda_{j-1}}{h_M},  \frac{\lambda_{j}}{h_M}\right] \right\} \geq 1 \right\},
    \end{equation*}
    that is valid and dominates \eqref{eq:ghom}.
\end{lemma}

\begin{proof}
    It is simple to see that $f$ is decreasing and upper semicontinuous. In addition,
    \[
    \int_0^1 f(p) \d p = \sum_{j=1}^M \frac{1}{\lambda_j}\left(\frac{\lambda_j - \lambda_{j-1}}{h_M} \right) = \frac{1}{h_M} \sum_{j=1}^M \frac{\lambda_j - \lambda_{j-1}}{\lambda_j} = 1.
    \]
    This implies that $F_{\mathrm{GHom}}$ is admissible.
To prove the last part, we see that for any $\mathbf{p} \in (0,1]^K$ and $\alpha \in (0, 1)$ we have that $F_\mathrm{GHom}(\mathbf{p}) \leq \alpha$, if and only if,
    \[
    \bigwedge_{k=1}^M \frac{1}{\lambda_k}p_{(\lceil \lambda_k K \rceil)} \leq \frac{\alpha}{h_M}.
    \]
    This implies that, for some $m \in \{1, \dots, K\}$, we have that 
    \[
    \sum_{i=1}^K \frac{1}{K}\ind\left\{p_i \leq \alpha \frac{\lambda_m}{h_M} \right\} \geq \lambda_m.
    \]
    We now define this chain of inequalities,
    \[
    \begin{split}
        1 &\leq \sum_{i=1}^K \frac{1}{K\lambda_m} \ind\left\{p_i \leq \alpha\frac{\lambda_m}{h_M} \right\} = \frac{1}{K} \sum_{i=1}^K \frac{1}{\lambda_m} \sum_{j=1}^m \ind\left\{\frac{p_i}{\alpha} \in \left(\frac{\lambda_{j-1}}{h_M}, \frac{\lambda_j}{h_M} \right]\right\}\\
        & \leq \frac{1}{K} \sum_{i=1}^K  \sum_{j=1}^m \frac{1}{\lambda_j} \ind\left\{\frac{p_i}{\alpha} \in \left(\frac{\lambda_{j-1}}{h_M}, \frac{\lambda_j}{h_M} \right]\right\} \leq  \frac{1}{K} \sum_{i=1}^K  \sum_{j=1}^M \frac{1}{\lambda_j} \ind\left\{\frac{p_i}{\alpha} \in \left(\frac{\lambda_{j-1}}{h_M}, \frac{\lambda_j}{h_M} \right]\right\}.
    \end{split}
    \]
\end{proof}    

It is possible to see that \eqref{eq:ghom} is a special case of the R\"uger combination when $M=1$, while it coincides with the Hommel combination rule when $\lambda = (0, 1/K, \dots, (K-1)/K, 1)$.

\section{Improving generalized mean}\label{sec:gen_mean}
In this section, we discuss the generalized mean combination rule, for $r \in \mathbb{R}\setminus\{0\}$. This combination rule, introduced in~\cite{vovk2020}, is quite broad and contains some important cases well known in the literature. In particular, if $r=1$ it reduces to the sample average (Section~\ref{sec:avg}), while if $r=-1$ it coincides with the harmonic mean described in Section~\ref{sec:harm_avg}. We introduce a lemma characterizing the calibrator used in the context of the generalized mean combination rule.
\begin{lemma}\label{lemma:calib_gen_avg}
    Let $r \in \mathbb{R}\setminus \{0\}$ and $f:[0,\infty) \to [0, \infty]$ be given by
    \begin{equation}\label{eq:calib_gen_avg}
        f(p) = \min\left\{\frac{r(1-p^r)}{T_{r,K}}, K \right\} \ind\{ p  \in [0,1] \}, 
    \end{equation}
    where $T_{r,K} > 0$ is any constant, possibly dependent on $K$, such that $\int_0^1 f(p) \d p \leq 1$. Then $f$ is a calibrator.
\end{lemma}
\begin{proof}
    Since $f(p)$ is continuously decreasing in $T_{r,K}$, it can be verified that there exists $T_{r,K} > 0$ such that $\int_0^1 f(p) \d p = 1$. Moreover, $f$ is decreasing and non-negative which completes the claim.
\end{proof}

It is simple to see that for $r>0$, we have that $T_{r,K}=r^2/(r+1)$ satisfies $\int_0^1 f(p) \d p \leq 1$. From previous sections, we obtain $T_{1,K} = 1/2$ while a more complex result appears for $T_{-1,K}$. We now see how the calibrator defined in~\eqref{eq:calib_gen_avg} is related to the rule defined in Section \ref{sec:general_res}. First, we define the function $F'_\mathrm{M_r}$ as
\begin{equation}\label{eq:FMr2}
    F'_\mathrm{M_r}(\mathbf{p}) = \inf\left\{\alpha \in (0,1): \frac{1}{K} \sum_{k=1}^K \frac{r(1-(p_k/\alpha)^r)}{T_{r,K}} \geq 1 \right\} = \frac{M_\mathrm{r}(\mathbf{p})}{(1 - T_{r,K}/r)^{1/r}},
\end{equation}
where $M_\textrm{r}(\mathbf{p})$ is the $r$-generalized mean of $\mathbf{p}$, defined by $M_\textrm{r}(\mathbf{p})=(\sum_{k=1}^K p_k^r/K)^{1/r}$. In particular, $F'_\mathrm{M_r}(\mathbf{p})$ coincides with the generalized mean combination rule where $a_{r,K}=(1-T_{r,K}/r)^{-1/r}$. In addition, if $r>0$ then $F'_\mathrm{M_r}(\mathbf{p})=(r+1)^{1/r}{M_r}(\mathbf{p})$ which coincides with the asymptotically precise merging function studied by~\cite{vovk2020}. It is possible to prove that $F'_\mathrm{M_r}(\mathbf{p}) \leq F_\mathrm{M_r}(\mathbf{p})$, where 
\begin{equation}\label{eq:FMr}
    F_\mathrm{M_r}(\mathbf{p}) = \inf\left\{\alpha \in (0,1): \frac{1}{K} \sum_{k=1}^K\left( \frac{r(1-(p_k/\alpha)^r)}{T_{r,K}}\right)_+ \geq 1 \right\},
\end{equation}
which is the p-merging function induced by the calibrator defined in~\eqref{eq:calib_gen_avg}. In particular, according to Theorem~\ref{th:vovkwang} we have that $F_\mathrm{M_r}(\mathbf{p})$ is a p-merging function.

\subsection{Exchangeable generalized mean}
We now define the function $F_{\mathrm{EM_r}}$ in the following way:
\begin{equation}\label{eq:FEMr}
    F_{\mathrm{EM_r}}(\mathbf{p}) = \inf\left\{\alpha \in (0,1): \bigvee_{\ell \le K} \frac{1}{\ell} \sum_{i=1}^\ell \left(\frac{r(1-(p_i/\alpha)^r)}{T_{r,K}} \right)_+ \geq 1\right\}.
\end{equation}
If we define $F'_{\mathrm{EM_r}}$ by
\begin{equation}\label{eq:FEMr2}
    F'_{\mathrm{EM_r}}(\mathbf{p}) = \inf\left\{\alpha \in (0,1): \bigvee_{\ell\le K} \frac{1}{\ell} \sum_{i=1}^\ell \left(\frac{r(1-(p_i/\alpha)^r)}{T_{r,K}} \right) \geq 1\right\},
\end{equation}
then it holds that $F_{\mathrm{EM_r}} \leq F'_{\mathrm{EM_r}}$.
\begin{theorem}
    Let $r\in \mathbb{R}\setminus \{0\}$, then the function $F'_{\mathrm{EM_r}}$ defined in~\eqref{eq:FEMr2} equals 
    \[
    \left(1-\frac{T_{r,K}}{r} \right)^{-1/r} \left( \bigwedge_{m=1}^K M_r(\mathbf{p}_m) \right),
    \]
    is an ex-p-merging function, and it strictly dominates the function $F_{\mathrm{M_r}}$ in~\eqref{eq:FMr}. However, it is strictly dominated by the function $F_{\mathrm{EM_r}}$ defined in~\eqref{eq:FEMr} that is also an ex-p-merging function. 
\end{theorem}
\begin{proof}
    According to Theorem~\ref{th:calibr_to_F_exc} and Lemma~\ref{lemma:calib_gen_avg}, it follows that the function $F_{\mathrm{EM_r}}$ is an ex-p-merging function. In addition, fix any $\alpha \in (0,1)$ and $\mathbf{p} \in (0,1]^{K}$, then $F'_{\mathrm{EM_r}} \leq \alpha$ if and only if  
    \[
    \frac{1}{\ell} \sum_{i=1}^\ell \frac{r(1-(p_i/\alpha)^r)}{T_{r,K}} \geq 1 \mathrm{~for~some~}\ell \leq K \implies \left(1-\frac{T_{r,K}}{r}\right)^{-1/r} \left(\frac{1}{\ell} \sum_{i=1}^\ell p_i^r\right)^{1/r} \leq \alpha \mathrm{~for~some~}\ell \leq K.
    \]
    Taking a minimum over $\ell$ yields the desired formula.
\end{proof}

\subsection{Randomized generalized mean}
According to the the previous sections, we define the randomized p-merging function $F_{\mathrm{UM_r}}$ as follows
\begin{equation}\label{eq:FUMr}
    F_\mathrm{UM_r}(\mathbf{p}, u) = \inf\left\{\alpha \in (0,1): \frac{1}{K} \sum_{k=1}^K\left( \frac{r(1-(p_k/\alpha)^r)}{T_{r,K}}\right)_+ \geq u \right\}.
\end{equation}
The function $F_{\mathrm{UM_r}} \leq F'_{\mathrm{UM_r}}$ where the function $F'_{\mathrm{UM_r}}$ is defined by
\begin{equation}\label{eq:FUMr2}
    F'_\mathrm{UM_r}(\mathbf{p}, u) = \inf\left\{\alpha \in (0,1): \frac{1}{K} \sum_{k=1}^K\left( \frac{r(1-(p_k/\alpha)^r)}{T_{r,K}}\right) \geq u \right\} ,
\end{equation}
and can be considered as the randomized version of~\eqref{eq:FMr2}.
\begin{theorem}
    Let $r\in \mathbb{R}\setminus \{0\}$, then the function $F'_{\mathrm{EM_r}}$ defined in~\eqref{eq:FUMr2} equals 
    \[
    \frac{M_r(\mathbf{p})}{(1-uT_{r,K}/r)^{1/r}},
    \]
    is a randomized p-merging function, and it strictly dominates the function $F_{\mathrm{M_r}}$ in~\eqref{eq:FMr}. However, it is strictly dominated by the function $F_{\mathrm{UM_r}}$ defined in~\eqref{eq:FUMr} that is also a randomized p-merging function. 
\end{theorem}
\begin{proof}
According to Corollary~\ref{co:1} and Lemma~\ref{lemma:calib_gen_avg}, it follows that $F_{\mathrm{UM_r}}$ is a valid randomized p-merging function. In addition, fix any $\alpha \in (0,1)$ and $\mathbf{p} \in (0,1]^K$, then $F'_{\mathrm{UM_r}} \leq \alpha$ if and only if 
\[
\frac{1}{K} \sum_{k=1}^K\left( \frac{r(1-(p_k/\alpha)^r)}{T_{r,K}}\right) \geq u \implies \frac{M_r(\mathbf{p})}{(1-uT_{r,K}/r)^{1/r}} \leq \alpha.
\]
\end{proof}

\section{General algorithm}\label{subsec:algo}
One general algorithm to compute the ex-p-merging function defined in~\eqref{eq:f_exch} induced by a calibrator $f$ is proposed in Algorithm~\ref{alg:alpha_search}. The algorithm employs the bisection method and it consistently yields a p-value exceeding that of the induced ex-p-merging function by at most $2^{-B}$.
\begin{algorithm}[H]
\caption{Ex-p-merging function}
\label{alg:alpha_search}
\begin{algorithmic}
    \Require A calibrator $f$, $B \in \mathbb{N}$, and a sequence of p-values $(p_1, \dots, p_K)$
\State $L := 0$ and $U := 1$
\For{$m = 1, \dots, B$}
    \State $\alpha := (L+U)/2$
    \If{$\bigvee_{k=1}^K \left(\frac{1}{k} \sum_{i=1}^k f\left( \frac{p_i}{\alpha} \right) \right) \geq 1$}
        \State $U := \alpha$
    \Else
        \State $L := \alpha$
    \EndIf
\EndFor
\Ensure $U$
\end{algorithmic}
\end{algorithm}

\section{Additional simulation results}\label{sec:add_sims}
In this section we report some additional simulation results. We first compare the merging rules introduced in Section~\ref{sec:ruger} and in Section~\ref{sec:avg}, specifically the rules: ``twice the median'' and ``twice the average''. In particular, $K=100$ p-values are generated as in \eqref{eq:sim_pvals}, and two different values of $\rho$ are chosen, 0.9 and 0.1, respectively. The parameter $k$ for the randomized R\"uger combination rule is set to $K/2$ and $\mu$ varies in the interval $[0,3]$. The results are computed by averaging the outcomes obtained in  $10,000$ replications.

\begin{figure}[t!]
    \centering
    \includegraphics[width = 1\textwidth]{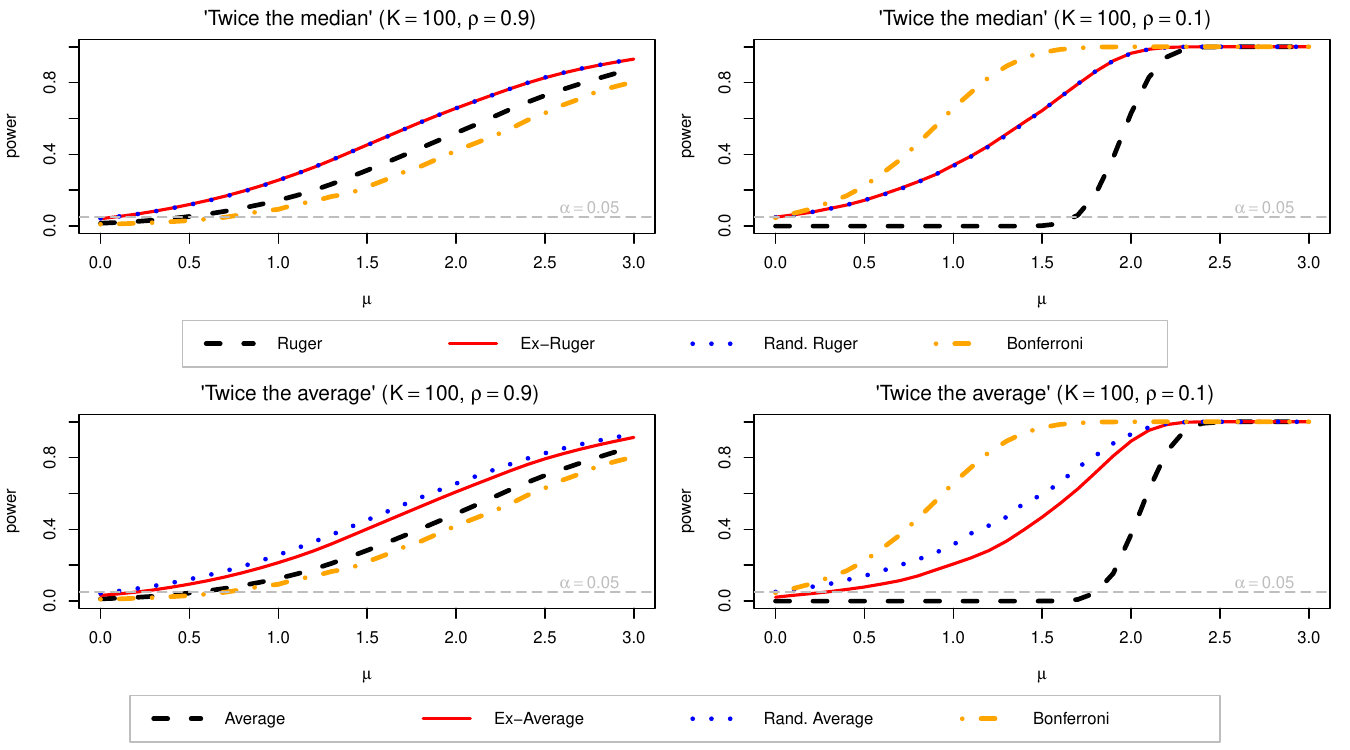}
    \caption{Combination of p-values using different rules. Every subplot illustrates power against $\mu$. The left endpoint of $\mu = 0$ actually represents the empirical type I error, which is controlled at the nominal level $\alpha = 0.05$ for all methods proposed. The first column has $\rho=0.9$, while the second column has $\rho=0.1$ --- as expected, the Bonferroni correction is more powerful near independence, but is less powerful under strong dependence. Further, our exchangeable and randomized improvements offer sizeable increases in power over the original variants in all settings.}
    \label{fig:sims}
\end{figure}

In Figure~\ref{fig:sims}, we can see that the type I error is controlled at the nominal level $0.05$ for all proposed methods. In the case of the R\"uger combination, in both dependence scenarios, the power of the combinations obtained by exploiting exchangeability or employing external randomization is highly similar. 
In the case of the combination based on the arithmetic mean, it appears that the rules obtained using randomization exhibit greater power than those based on exchangeability (and, naturally, than the original rules). In general, across all observed scenarios, the new rules demonstrate a quite significant improvement in terms of power. 

In addition, we compare all the randomized combination rules reported in the paper. We omit the randomized functions that are dominated by other randomized p-merging functions. As before, $K=100$ p-values are generated as in \eqref{eq:sim_pvals}, and $\rho=\{0.1, 0.9\}$. The parameter $k$ for the randomized R\"uger combination rule is set to $K/2$ and $\mu$ varies in the interval $[0,3]$. The results are obtained by repeating the procedure $10,000$ times and reporting the average.
\begin{figure}
    \centering
    \includegraphics[width=\textwidth]{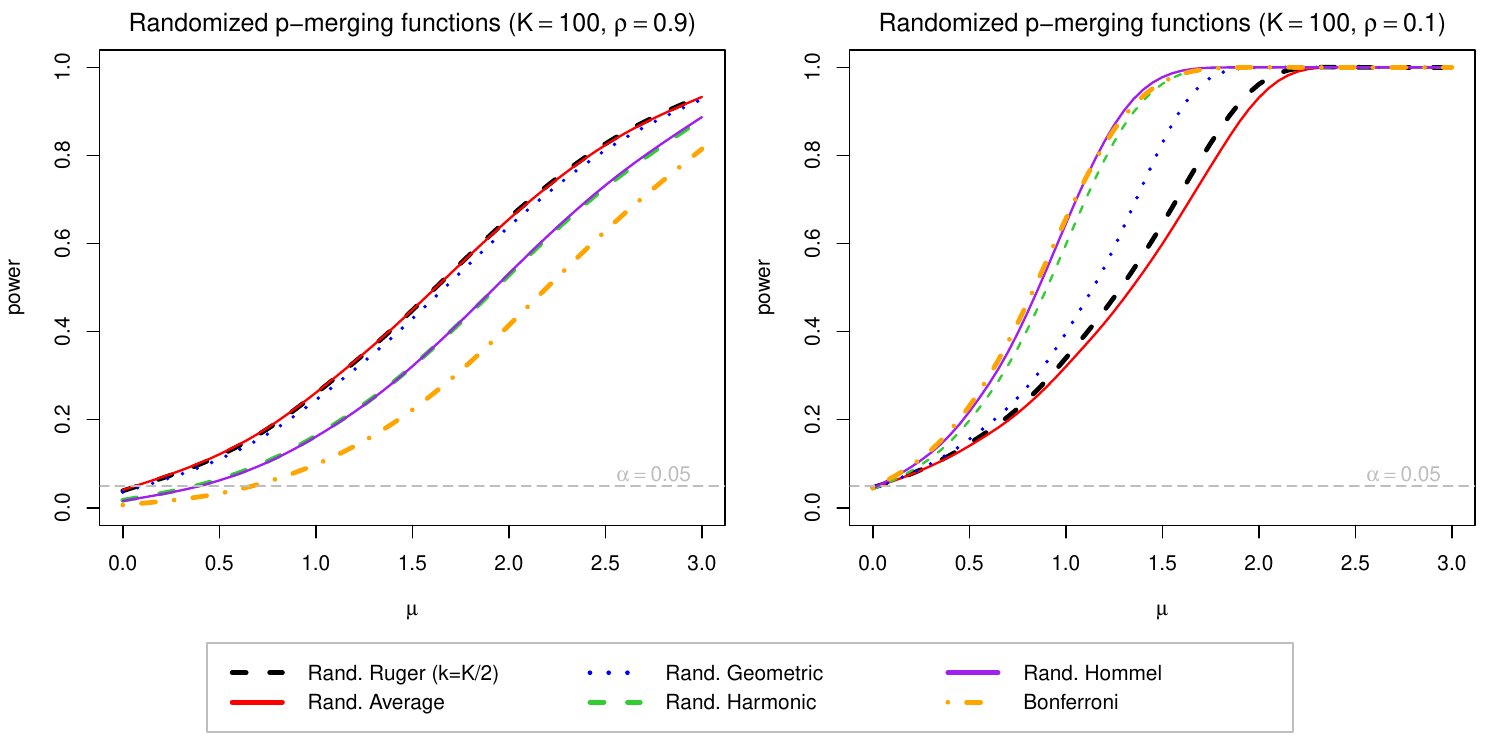}
    \caption{Combination of p-values using different randomized p-merging functions. The order of the performance of the different ex-p-merging functions is almost the opposite in the two situations.}
    \label{fig:rand_funs}
\end{figure}

The results of the randomized versions of ``twice" the median and (improved) ``twice" the average are similar in both scenario. In addition, $F_{\mathrm{UHom}}$ and $F_{\mathrm{UH}}$ are quite similar in terms of power (Figure~\ref{fig:rand_funs}). From Figure~\ref{fig:rand_funs}, we can observe an opposite behavior of the functions in the case where $\rho = 0.9$ or $\rho = 0.1$. In general, the most powerful functions in the left graph tend to be the least powerful in the right graph, and vice versa.

At the end, we examine the 3 different randomized combination rules defined in Section~\ref{sec:avg}, respectively, $F_{\mathrm{UA}}, F'_{\mathrm{UA}}, F^*_{\mathrm{UA}}$. In particular, we recall that, $F_{\mathrm{UA}}$   dominates $F'_\mathrm{{UA}}$, while the combination $F^*_{\mathrm{UA}}$ (introduced in \citet[Appendix B.2]{wang2024testing}) neither dominates nor is dominated by either of the two.
\begin{figure}
    \centering
    \includegraphics[width=\textwidth]{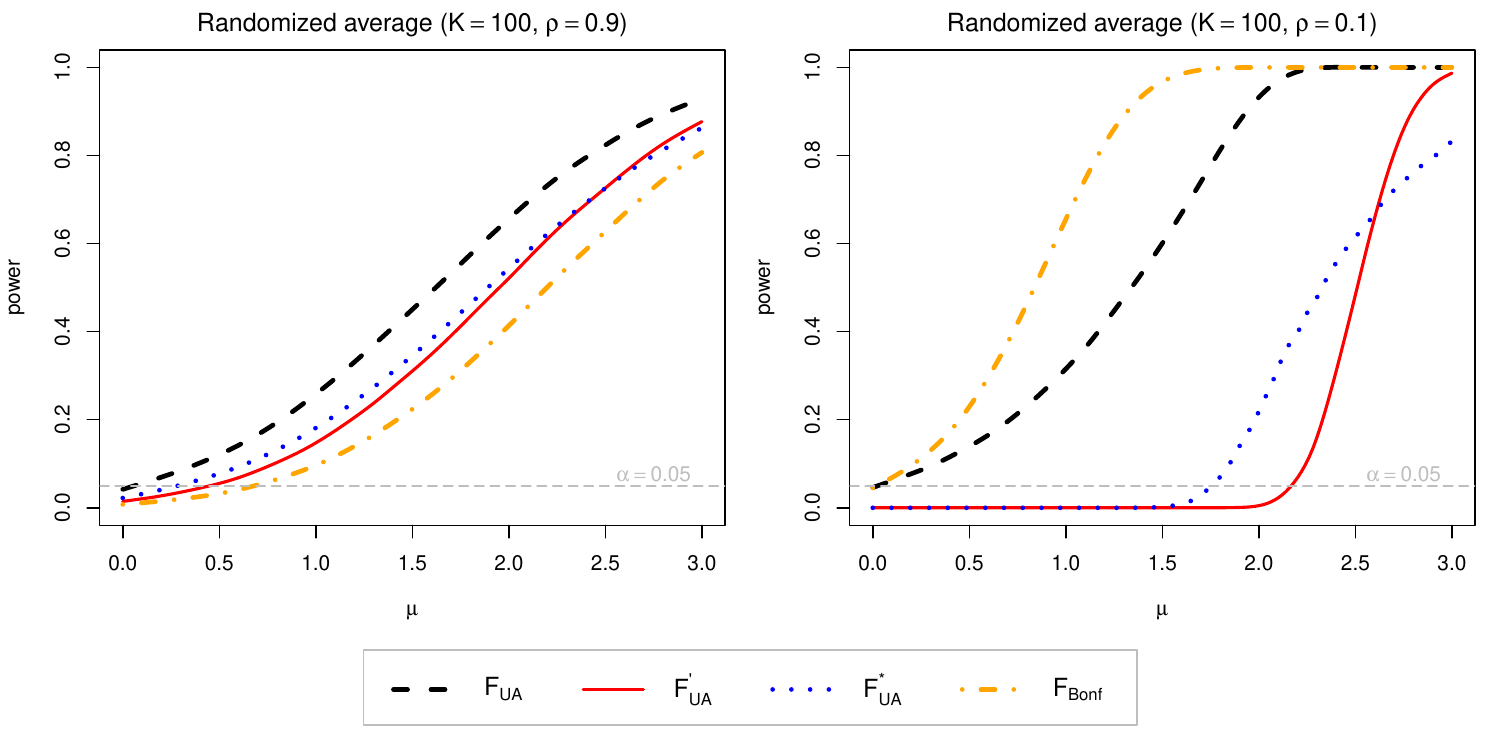}
    \caption{Combination of p-values using different randomized combination rules based on the average of p-values (and the Bonferroni method, for comparison). $F^*_{\mathrm{UA}}$ is more powerful than $F'_{\mathrm{UA}}$ only when $\mu \lesssim 2.5$.}
    \label{fig:rand_avg}
\end{figure}

The set of $K=100$ p-values is generated as in \eqref{eq:sim_pvals}, and two different values of $\rho$ are chosen, 0.9 and 0.1, respectively. The results are obtained by repeating the procedure $10,000$ times and reporting the average. From Figure~\ref{fig:rand_avg}, we can see that the power $F_{\mathrm{UA}}$ is always higher than the power of the other two combination rules. The function $F^*_{\mathrm{UA}}$ is more powerful than $F'_{\mathrm{UA}}$ only in the first part of the $y$-axis in both cases (for $\mu \leq 2.5$, indicatively), so there is no clear preference between the two. As expected, the Bonferroni correction is more powerful near independence ($\rho = 0.1$).

Before concluding, let us see what our procedures give for testing a global null. The goal of this last part is to explore a case where the order of the p-values in our proposed ex-p-merging functions is chosen in a data-driven manner. In the considered scenario, we will see how a particular statistics can be chosen to order the p-values with the aim of improving the statistical power under the alternative. In particular, this data-driven ordering does not alter the exchangeability under the null that is required in our ex-p-merging functions. Specifically, we investigate the issue of performing simultaneous one-sample t-tests using $n$ observations for each hypothesis. Let us suppose to have $K$ samples (one for each hypothesis) from a normal distribution,
$$
X_{ki} \sim \mathcal{N}(\mu_k, \sigma^2_k), \quad i=1,\dots,n,\quad \mu_k \in \mathbb{R}, \quad \sigma_k>0,
$$ 
where the observations $X_{ki}$, $k \in \{1,\dots,K\}$, $i\in\{1,\dots,n\}$, are mutually independent. Our goal is to test the hypothesis $H_k:\mu_k=0$ and to do so we use t-test. In addition, we define the global null as $H_0: \bigcap_{k=1}^K H_k$, that can be tested by merging the different p-values into a single p-value. The same problem has been studied, for example, in \cite{westfall2004} and \cite{ignatiadis2024}. Let us define
\[
\Bar{X}_k=\frac{1}{n} \sum_{i=1}^n X_{ki}, \quad \hat{\sigma}_k^2=\frac{1}{n-1} \sum_{i=1}^n (X_{ki} - \Bar{X}_k)^2, \quad T_k = \frac{\sqrt{n} \Bar{X}_k}{\hat{\sigma}_k^2},
\]
and p-values are $P_k:= 2G_{n-1}(-|T_k|)$, where $G_{n-1}$ is the t-distribution with $n-1$ degrees of freedom.

Our p-values will be ordered using the statistics $S_k^2 = \sum_{i=1}^n X_{ki}^2$, and this can be done since the proposed statistics are independent from $T_k$ (and so $P_k$). The key observation is that, under $H_0$ (i.e., when $\mu_k=0$ for all $k=1,\dots,K$), then the statistic $S_k^2$ is sufficient and complete for the inference on the parameter $\sigma_k^2$. On the other hand, the test statistic $T_k$ is constant in distribution with respect to $\sigma_k^2$. By Basu's Theorem \citep{basu1955}, $S_k^2$ and $T_k$ are independent and, in addition, the test $T_k$ and the statistics $S_k^2$ are independent among themselves since they are functions of independent random variables.

\begin{figure}
    \centering
    \includegraphics[width=\linewidth]{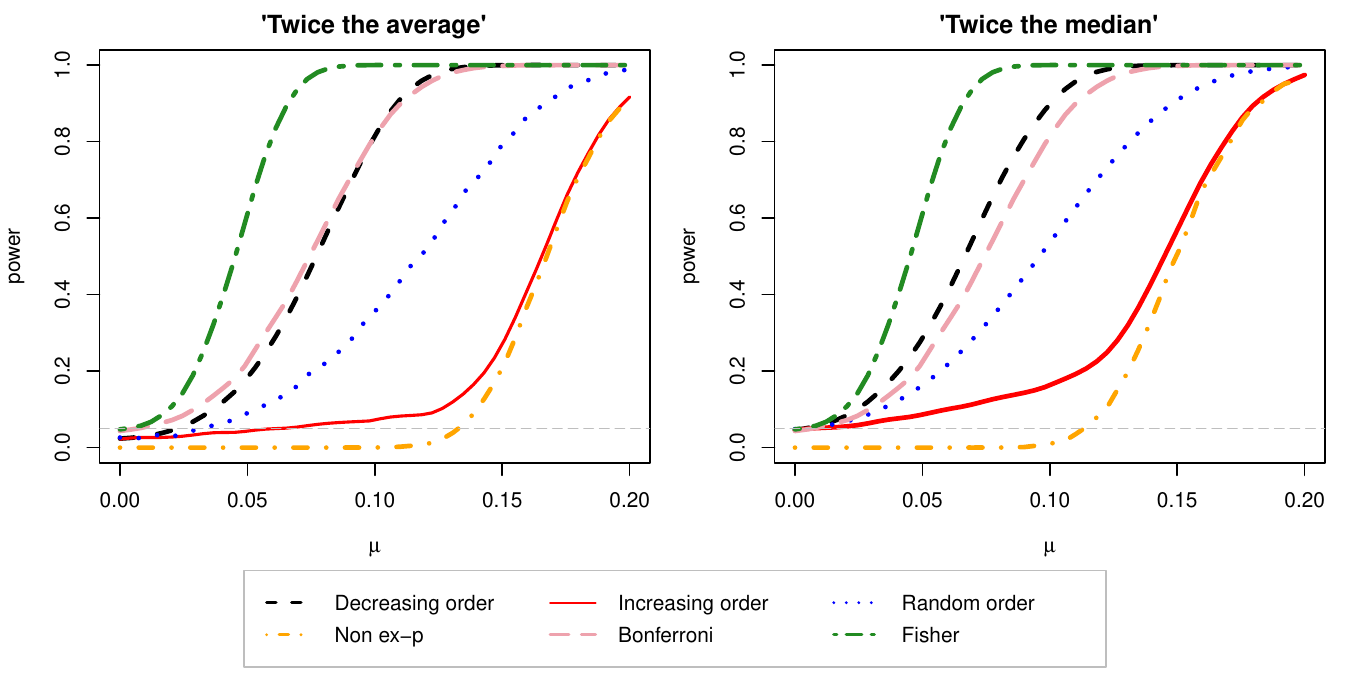}
    \caption{Combination of p-values for testing a global null. If p-values are ordered in decreasing order with respect to $S_k^2$ we have a power that is comparable with the Bonferroni rule. However, in all situations Fisher's combination is the most powerful.}
    \label{fig:global}
\end{figure}

We assume $\sigma_k=\sigma=1$ for all $k=1,\dots,K$, and $\mu_k=k\cdot\mu$, where $\mu$ is a parameter that varies in $[0,0.2]$. The parameter $\alpha$ is set to $0.05$, while $K=20$ and $n=10$. As can be seen in Figure \ref{fig:global}, in all situations the type I error is controlled at the level $\alpha$ under $H_0$ (i.e., when $\mu=0$). If p-values are ordered in decreasing order with respect to the statistics $S_k^2$ then we have a rule that has a power comparable (or higher) than the Bonferroni rule. It is expected that the descending order is more effective, as a large $S_k^2$ indicates that $H_k$ is likely false. The Fisher combination, which requires the strong assumption of independence, has the largest power.

\end{document}